\documentclass[AMA,STIX1COL]{WileyNJD-v2}
\usepackage{caption}
\usepackage{color}
\usepackage{placeins}
\usepackage{amsmath,amsfonts,amsthm,epsfig}
\usepackage{booktabs}
\usepackage{titlesec}
\usepackage{tikz}
\usepackage{pgfplots}

\pgfplotsset{compat=1.11,
    /pgfplots/ybar legend/.style={
    /pgfplots/legend image code/.code={%
       \draw[##1,/tikz/.cd,yshift=-0.25em]
        (0cm,0cm) rectangle (3pt,0.8em);},
   },
}

\begin{document}

\title{Towards Accuracy and Scalability: Combining Isogeometric Analysis with Deflation to Obtain Scalable Convergence for the Helmholtz Equation.}

\author[1]{V. Dwarka}
\author[1]{R. Tielen}
\author[1]{M. M\"oller}
\author[1]{C. Vuik}

\authormark{V. Dwarka,  R. Tielen, M. M\"oller and C. Vuik}

\address[1]{\orgdiv{Delft University of Technology}, \orgname{Numerical Analysis}, \orgaddress{\state{Delft}, \country{the Netherlands}}}
\presentaddress{van Mourik Broekmanweg 6, 2628 XE, Delft, the Netherlands}
 
\abstract[Abstract]{Finding fast yet accurate numerical solutions to the Helmholtz equation remains a challenging task. The pollution error (i.e. the discrepancy between the numerical and analytical wave number $k$) requires the mesh resolution to be kept fine enough to obtain accurate solutions. A recent study showed that the use of Isogeometric Analysis (IgA) for the spatial discretization significantly reduces the pollution error. 

However, solving the resulting linear systems by means of a direct solver remains computationally expensive when large wave numbers or multiple dimensions are considered. An alternative lies in the use of (preconditioned) Krylov subspace methods. Recently, the use of the exact Complex Shifted Laplacian Preconditioner (CSLP) with a small complex shift has shown to lead to wave number independent convergence while obtaining more accurate numerical solutions using IgA. 

In this paper, we propose the use of deflation techniques combined with an approximated inverse of the CSLP using a geometric multigrid method. Numerical results obtained for both one- and two-dimensional model problems, including constant and non-constant wave numbers, show scalable convergence with respect to the wave number and approximation order $p$ of the spatial discretization. Furthermore, when $kh$ is kept constant, the proposed approach leads to a significant reduction of the computational time compared to the use of the exact inverse of the CSLP with a small shift.}

\keywords{Helmholtz, pollution, numerical dispersion, Isogeometric Analysis, FEM, GMRES, Deflation}

%\maketitle 

\section{Introduction}
The Helmholtz equation has been widely studied in various fields of physics ranging from biomedical physics to geo- and nuclear physics. The electromagnetic scattering problem thus finds many applications in engineering practices. Many efforts have been made to find fast yet accurate numerical solutions to the Helmholtz problem. The latter remains a challenging topic in research due to the pollution error and the resulting linear system having undesirable properties. In particular, the pollution error results from a discrepancy between the analytical and numerical wave number \cite{ihlenburg1995dispersion,ihlenburg1997finite,ihlenburg1997solution}. Consequently, the mesh resolution has to be kept fine enough to obtain accurate numerical solutions. If we let $k$ denote the wave number, $N_{\rm dof}$ the number of degrees of freedom in one-dimension and $p$ the order of a finite difference or standard finite element scheme, then 
\begin{align*}
    N_{\rm dof} = Ck^{\left( \frac{p+1}{p} \right)},
\end{align*}
where $C$ is a constant that only depends on the accuracy achieved \cite{turkel2013compact}. In practice, this has led to the rule of thumb $kh \approx \frac{2\pi}{10}$, where $10$ denotes the number of degrees of freedom per wavelength and $h$ the mesh width. However, the resulting numerical solution still suffers from pollution, unless the resolution is kept at $C(k^{p+1} h^p) \leq 1$, for a general $p-$th order scheme. While this minimizes the pollution error, the resulting linear systems are too large for direct solution methods. This exacerbates in higher-dimensions, which opens the door to the use of iterative solution methods. Due to the resulting linear systems being indefinite and non-Hermitian, Krylov subspace or Induced Dimension Reduction methods are necessary. In fact, even using standard multigrid as a stand-alone solver diverges for the Helmholtz equation \cite{ernst2012difficult,ernst2013}. Moreover, for Krylov subspace methods, the number of iterations until convergence grows with the wave number $k$. Thus, the difficulty in solving Helmholtz-type problems can be reduced to optimizing the trade-off between having accurate numerical solutions, while using a scalable solver. 

One potential way to mitigate this problem is to adopt Isogeometric Analysis (IgA) \cite{hughes2005} as a discretization technique. IgA can be considered as the natural extension of the finite element method (FEM) to higher-order B-splines and has become widely accepted as a viable alternative to standard FEM. The use of high-order B-splines or Non-Uniform Rational B-splines (NURBS) enables a highly accurate representation of complex geometries and bridges the gap between computer-aided design (CAD) and computer-aided engineering (CAE) tools. Furthermore, a higher accuracy per degree of freedom can be achieved compared to standard FEM \cite{HUghes2007}. A new branch of studies has demonstrated that IgA furthermore helps to control the pollution error while keeping the size of the resulting linear system moderate \cite{buffa2010isogeometric,buffa2014isogeometric,wu2015isogeometric,coox2016performance,drzisga2020surrogate}. In \cite{mederos2020isogeometric}, the authors investigated the obtained accuracy for several Helmholtz-type problems using a non-constant wave number and documented increased accuracy. Thus, while the use of IgA for Helmholtz-type problems becomes more established, the process of solving the underlying discretized systems remained fairly untouched. Until recently, a study by Diwan et al. \cite{diwan2020iterative} covered this for the Helmholtz equation and researched the use of IgA together with an iterative solver. There, the resulting linear systems are solved using the Generalized Minimum Residual Krylov method (GMRES) preconditioned with the Complex Shifted Laplacian Preconditioner (CSLP) using a small complex shift. The results show wave number independent convergence of the iterative solver and, at the same time, higher accuracy of the numerical solution.

The well-known CSLP has been the industry standard for many years \cite{erlangga2006novel}. While this has accelerated the convergence dramatically, the number of iterations increases with the wave number $k$, which is why in order to obtain wave number independent convergence, the complex shift has to be kept at $\mathcal{O}(k^{-1})$ \cite{gander2015applying}. One drawback of keeping the shift very small is that the resulting preconditioner starts resembling the original matrix and exact inversion puts a heavy tax on the computational resources. Therefore, a few multigrid cycles are often used to approximate the inverse of the CSLP, which amounts to $\mathcal{O}(N)$ FLOPs \cite{erlangga2006novel}. However, in order to prevent multigrid from diverging, the complex shift has to be kept as large as possible $\mathcal{O}(1)$ \cite{cocquet2017large}. 

As a consequence, recent developments have led to a broad range of preconditioners such as domain decomposition based preconditioners \cite{gander2013domain,gander2014optimized,graham2017recent,bonazzoli2017two,bonazzoli2019domain,bootland2019dirichlet,graham2020domain}, sweeping preconditioners \cite{engquist2011sweeping,liu2016recursive,stolk2017improved,gander2019class,taus2020sweeps} and (multilevel) deflation based preconditioners \cite{sheikh2016accelerating,erlangga2017multilevel,dwarka2020scalable}. One of these new preconditioners is the Adapted Deflation Preconditioner (ADP), which uses higher-order Bezier curves to construct the deflation space. For finite difference discretizations, the preconditioner has shown to be simple yet competitive to the small-shift and exact inversion of CSLP in terms of wave number independent convergence and computational complexity for large wave numbers $k$. In essence, the deflation preconditioner projects the near-zero eigenvalues of the CSLP-preconditioned system onto zero. These near-zero eigenvalues are known to interfere with fast convergence of the Krylov subspace solver. 
% -paper superior wrt pollution.
% -first step with respect to iterative solution methods and preconditioning
% -in this paper we extend this research direction and focus on scalable convergence and efficiency in terms of CPU timings. 

Consequently, our aim in this paper is to extend the research direction set out in \cite{diwan2020iterative,mederos2020isogeometric}, by combining state-of-the-art iterative solvers with IgA discretization techniques to obtain both accurate and computationally efficient numerical solutions. In particular, we propose the use of deflation techniques combined with an approximated inverse of the CSLP using multigrid to obtain scalable and faster convergence with respect to the wave number $k$ and the order $p$. 
We study one- and two-dimensional model problems using IgA discretizations containing both a constant wave number $k$ and a variable wave number $k(x,y)$. In the latter case, we focus on the performance of the solver in the presence of sharp discontinuities in the wave number and the underlying solution. 
For the two-dimensional model problems, we report the number of iterations and the CPU-timings to show that the use of deflation combined with a multigrid-approximated CSLP allows for tremendous gain in computational efficiency while keeping scalable convergence in terms of the number of iterations. The method outperforms the exact inversion of the CSLP with a small complex shift in terms of number of iterations and CPU-timings when a large constant or non-constant wave number is used. 

The paper is organized as follows. We start with the variational formulation of the Helmholtz equation and the model problem definitions in section \ref{probdef}. In section \ref{defgmres} we discuss the deflation preconditioning technique for the Krylov subspace method. Here we introduce the use of higher-order Bezier curves as a basis for the deflation space. We then proceed by performing a spectral analysis of the preconditioned systems and various numerical experiments in section \ref{numres} in order to determine the convergence behavior. We provide CPU-timings in order to assess the computational time complexity. We conclude our results in section \ref{concl}. 

% \begin{enumerate}
%     \item Introduction to the Helmholtz problem. - CHECK
%     \item Challenges: balancing between accuracy and fast convergence. Especially with respect to scalability. - CHECK
%     \item Brief history of pollution error. - CHECK
%     \item Reference to IgA and Helmholtz paper: short description of IgA literature and main motivation for using it in this context. 
%     \item Combining IgA with the newest/more practical sequential and iterative Helmholtz solver. - CHECK
%     \item Brief history of iterative Helmholtz solver: CSLP, deflation, domain decomposition etc. 
%     \item Emphasize bringing both topics together; accuracy and scalability. Give explicit definition of scalability in this context. - CHECK
%     \item End with paper overview. - CHECK
% \end{enumerate}

\section{Problem Definition}\label{probdef}
In order to assess the quality of the proposed solution method, we start by defining a variety of one- and two-dimensional model problems. In particular, we consider model problems involving both constant and non-constant wave numbers.  
Then, we proceed by presenting the variational formulation and B-spline discretization using the generalization of our two-dimensional model problem as an example.  

\subsection{One-dimensional model problems}
\subsubsection{MP 1-A}
The first one-dimensional model problem, MP 1-A, is given below
\begin{eqnarray} \label{MP1A}
-\frac{d^2 u(x)}{d x^2} -k^2 \, u &=&  0, \quad x \in  \Omega  = [0,1],  \\ 
 u(x)&=& 1, \quad x = 0, \nonumber \\
u'(x) - iku(x)&=& 0, \quad  x = 1. \nonumber 
\end{eqnarray}
Here, homogeneous Dirichlet and Sommerfeld boundary conditions are applied on the left and right boundary, respectively. The exact solution for MP1-A is given by $u(x) = e^{ikx}$. Model problem MP 1-A will be adopted to investigate the pollution error for various values of the approximation order $p$ of the B-spline basis functions. It will also be used to perform a convergence factor study in  order to check the robustness of the solver. \\
\\
\subsubsection{MP 1-B}
Model problem MP1-B involves an inhomogeneous source term. Furthermore, Dirichlet boundary conditions are applied on both boundaries, resulting in the following model problem 
\begin{eqnarray} \label{MP1B}
-\frac{d^2 u}{d x^2} -k^2 \, u &=&  \delta(x - x'), \quad x \in  \Omega  = [0,1], \\
 u(x) &=& 0, \qquad \qquad \ x=0, \nonumber \\
 u(x) &=& 0, \qquad \qquad \ x=1. \nonumber \\
 \end{eqnarray}
The analytic solution of MP1-B is based on the Green's function of this model problem and is given by 
\begin{eqnarray}
u(x,x') &=& 2 \sum\limits_{j=1}^{\infty}
\frac{ \sin\left( j \pi x \right)\sin\left( j \pi x' \right) }{j^2 \pi^2 - k^2}, \quad x \in  \Omega  = [0,1], \nonumber \\ 
k^2 &\neq& {j^2 \pi^2}, \hspace{3.3cm} j = 1,2,3,\hdots. \nonumber
\end{eqnarray}
Note that, for $k^2 = {j^2 \pi^2}$, the eigenfunction expansion would become defective as this would imply resonance and unbounded oscillations in the absence of dissipation. Therefore, we explicitly impose the extra condition $k^2 \neq {j^2 \pi^2}$ asserting that our Green's function exists. \\
By imposing Dirichlet boundary conditions, the resulting system matrix exhibits the most unfavorable distribution of the eigenvalues \cite{van2007spectral}. Note that the inclusion of Sommerfeld radiation conditions already slightly shifts the eigenvalues away from the origin due to the natural occurring damping. 

\subsection{Two-dimensional Model Problems}
\subsubsection{MP 2-A}
In two dimensions, we consider as MP 2-A the natural extension of MP 1-B to two dimensions:
 \begin{eqnarray} \label{MP2A}
 -\Delta u(x,y) - k^2u(x,y) &=& \delta{(x - \frac{1}{2},y - \frac{1}{2})}, \, \quad (x,y) \in \Omega = [0,1]^2,  \\
 u(x,y) &=& 0, \, \hspace{2.42cm} (x,y) \in \partial{\Omega_D},  \nonumber \\
 \left( \frac{\partial{}}{\partial{{\mathbf{n}}}} - ik \right)  u({x,y})  &=& {0}, \hspace{2.45cm} {(x,y)} \in \partial{\Omega_R}.
 \end{eqnarray}
 
 Again, the analytic solution is given by the Green's function:
 \begin{eqnarray} 
 u(x,y,x',y') &=& 4\sum\limits_{i=1}^{\infty} \sum\limits_{j=1}^{\infty} \frac{ \sin\left( i \pi x \right)\sin\left( i \pi x' \right) \sin\left( j \pi y \right)\sin\left( j \pi y' \right)}{i^2\pi^2 + j^2\pi^2  - k^2}, \quad (x,y) \in \Omega = [0,1]^2, \\
 k^2 &\neq& i^2\pi^2 + j^2\pi^2, \hspace{5.6cm} i,j = 1,2,3, \hdots. \nonumber
 \end{eqnarray}

\subsubsection{MP 2-B} 
As a final model problem, MP 2-B, we consider a non-constant wave number $k=k(x,y)$, an inhomogeneous source function and Dirichlet boundary conditions on the entire boundary $\partial \Omega$.
  \begin{eqnarray} \label{MP2B}
 -\Delta u(x,y) - k(x,y)^2u(x,y) &=& \delta{(x - \frac{1}{2},y - \frac{1}{2})}, \, \quad (x,y) \in \Omega = [0,1]^2,  \\
 u(x,y) &=& 0, \, \hspace{2.42cm} (x,y) \in \partial{\Omega}.  \nonumber 
 \end{eqnarray}

\begin{figure}[ht!]
    \centering
    \includegraphics[scale=0.45]{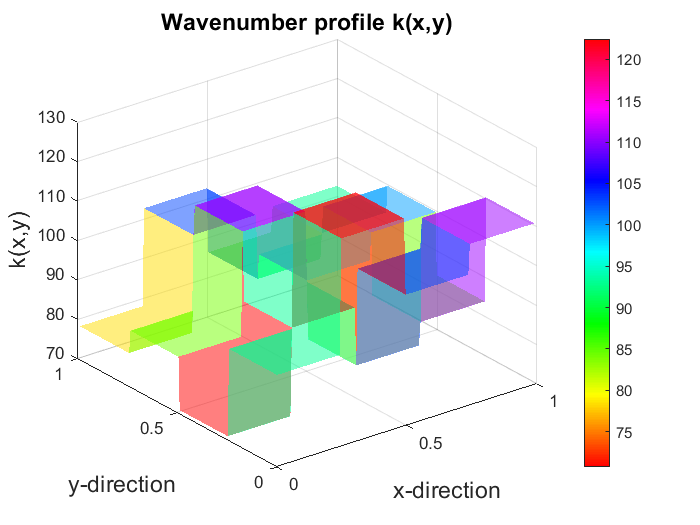}
    \caption{wave number distribution for $k(x,y)$. $k$ has been set to have a base value of 100. The figure shows the step-function to illustrate the variation profile of the wave number with respect to the $x-$ and $y-$direction.}
    \label{fig:kxy1}
\end{figure}

Here, $k(x,y)$ is chosen to be a two-dimensional step function consisting of $16$ different values. For a fixed value of $k$, the values vary between $\frac{1}{2}k$ and $\frac{3}{2}k$. Figure \ref{fig:kxy1} shows the considered field $k(x,y)$ for $k=100$. This model problem uses various horizontal layers in order to test the performance of the solver when a variable wave number $k(x,y)$ is used. This is particularly important to investigate as in certain cases for Helmholtz-type problems the underlying solver might diverge. This has been reported for domain decomposition based preconditioners using inexact factorizations \cite{gander2019class}.

\subsection{Variational Formulation}
\label{variational}

To illustrate the variational formulation, we consider the inhomogeneous Helmholtz equation in two dimensions adopting inhomogeneous Robin boundary conditions:
\begin{eqnarray} \label{eq:Helmholtz}
    \Delta u(x,y) - k^2 u(x,y) &=& f(x,y), \ \quad (x,y) \in \Omega \subset \mathbb{R}^2, \label{eq:testmp} \\
   \left (  \frac{\partial}{\partial \mathbf{n}} - ik \right ) u(x,y)  &=& g(x,y), \ \ \quad  (x,y) \in \partial \Omega. \label{eq:testmp2} 
\end{eqnarray} 

Here, $\Omega$ is a connected Lipschitz domain, $f \in L^2(\Omega)$, $g \in L^2(\partial \Omega)$ and $k > 0$ a constant wave number. Let us define $\mathcal{V} = H^1_0(\Omega)$ as the space of functions in the Sobolev space $H^1(\Omega)$ that vanish on the boundary $\partial \Omega$. The variational formulation of \eqref{eq:testmp} is obtained by multiplication with an test function $v \in \mathcal{V}$ and application of integration by parts

\begin{equation} \label{eq:bilinearform}
    a(u,v) = (f,v), \ \forall v \in \mathcal{V},
\end{equation}

where 

\begin{equation}
    a(u,v) = \int_{\Omega} \nabla u \cdot \overline{\nabla v} \ \text{d}\Omega + k^2 \int_{\Omega} u \overline{v} \ \text{d}\Omega - ik \int_{\partial \Omega} u \overline{v} \ \text{d} \Gamma  \hspace{1cm} (f,v) = \int_{\Omega} f\overline{v} \ \text{d}\Omega + \int_{\partial \Omega} g \overline{v} \ \text{d} \Gamma.
\end{equation}

A geometry function $\mathbf{F}$ is then defined to parameterize the physical domain $\Omega$ by describing an invertible mapping to connect the parameter domain $\Omega_0 = (0,1)^2$ with the physical domain $\Omega$.

\begin{equation}
    \mathbf{F} = \Omega_0 \rightarrow \Omega, \quad \mathbf{F}(\xi,\eta) = (x,y).
\end{equation}

The considered geometries throughout this paper can be described by a single geometry function $\mathbf{F}$, that is, the physical domain $\Omega$ is topologically equivalent to the unit square. In case of more complex geometries, a family of functions $\mathbf{F}^{(m)}$ ($m = 1,\ldots, K$) is defined and we refer to $\Omega$ as a multipatch geometry consisting of $m$ patches. 

\subsubsection{B-spline basis functions}
To discretize Equation \eqref{eq:testmp}, univariate B-spline basis functions are defined on the parameter domain $\Omega_0$ by an underlying knot vector $\Xi = \{ \xi_1, \xi_2, \ldots , \xi_{N+p}, \xi_{N+p+1} \}$. Here, $N$ denotes the number and $p$ the order of the B-spline basis functions. Based on this knot vector, the basis functions are defined recursively by the Cox-de Boor formula \cite{boor}, starting from the constant ones

\begin{eqnarray}
\phi_{j,0}(\xi) = \begin{cases} 1   \hspace{2.25 cm} \text{if} \hspace{0.2cm}\xi_j \leq \xi < \xi_{j+1}, \\  0 \hspace{2.25cm} \text{otherwise.}  \end{cases}
\end{eqnarray}

Higher-order B-spline basis functions of order $p>0$ are then defined recursively
\begin{eqnarray}
\phi_{j,p}(\xi) = \frac{\xi - \xi_j}{\xi_{j+p}-\xi_j} \phi_{j,p-1}(\xi) +  \frac{\xi_{j+p+1} - \xi}{\xi_{j+p+1}-\xi_{j+1}} \phi_{j+1,p-1}(\xi).
\end{eqnarray}

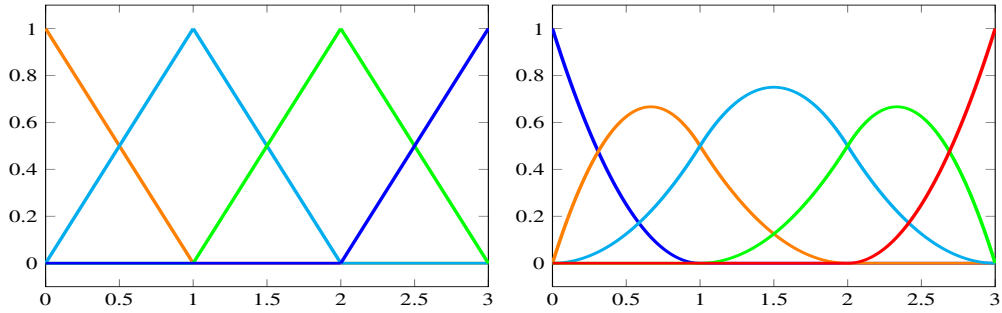
\begin{figure}[ht!]
\centering
\begin{tikzpicture}[xscale=0.85,yscale=0.65]
\begin{axis}[xmin=0,xmax=3,ymin=-0.1,ymax=1.1, samples=50]
%p=1
  \addplot[orange, ultra thick, domain=0:1] {1-x};
  \addplot[orange, ultra thick, domain=1:3] {0};
  \addplot[cyan, ultra thick, domain=0:1] {x};
  \addplot[cyan, ultra thick, domain=1:2] {2-x};
  \addplot[cyan, ultra thick, domain=2:3] {0};
  \addplot[green, ultra thick, domain=0:1] {0};
  \addplot[green, ultra thick, domain=1:2] {x-1};
  \addplot[green, ultra thick, domain=2:3] {3-x};
  \addplot[blue, ultra thick, domain=0:2] {0};
  \addplot[blue, ultra thick, domain=2:3] {x-2};
\end{axis}
\end{tikzpicture}
\begin{tikzpicture}[xscale=0.85,yscale=0.65]
\begin{axis}[xmin=0,xmax=3,ymin=-0.1,ymax=1.1, samples=50]
  \addplot[blue, ultra thick, domain=0:1] {(1-x)*(1-x)};
  \addplot[blue, ultra thick, domain=1:3] {0};
  \addplot[orange, ultra thick, domain=0:1] {x*(1-x)+0.5*(2-x)*x};
  \addplot[orange, ultra thick, domain=1:2] {0.5*(2-x)*(2-x)};
  \addplot[orange, ultra thick, domain=2:3] {0};
  \addplot[cyan, ultra thick, domain=0:1] {0.5*x*x+0.5*(3-x)*0};
  \addplot[cyan, ultra thick, domain=1:2] {0.5*x*(2-x)+0.5*(3-x)*(x-1)};
  \addplot[cyan, ultra thick, domain=2:3] {0.5*x*0+0.5*(3-x)*(3-x)};
  \addplot[green, ultra thick, domain=0:1] {0};
  \addplot[green, ultra thick, domain=1:2] {0.5*(x-1)*(x-1)};
  \addplot[green, ultra thick, domain=2:3] {0.5*(x-1)*(3-x)+(3-x)*(x-2)};
  \addplot[red, ultra thick, domain=0:2] {0};
  \addplot[red, ultra thick, domain=2:3] {(x-2)*(x-2)};
\end{axis}
\end{tikzpicture}
\caption{Linear and quadratic B-spline basis functions based on the knot vectors $\Xi_1 = \{ 0, 0, 1, 2, 3,3 \}$ and $\Xi_2 = \{ 0, 0, 0, 1, 2, 3,3,3 \}$, respectively.}
\label{fig:bspline}
\end{figure}

The resulting B-spline basis functions $\phi_{j,p}$ are non-zero on the interval $[\xi_j,\xi_{j+p+1})$ and possess the partition of unity property. Furthermore, the basis functions are $C^{p-m_j}$-continuous, where $m_j$ denotes the multiplicity of knot $\xi_j$. Throughout this paper, we consider a uniform knot vector with knot span size $h$, where the first and last knot are repeated $p+1$ times. As a consequence, the resulting B-spline basis functions are $C^{p-1}$ continuous and interpolatory at both end points. Figure \ref{fig:bspline} illustrates both linear and quadratic B-spline basis functions based on such a knot vector. 

For the two-dimensional case, the tensor product of univariate B-spline basis functions are adopted for the spatial discretization. Let $N_{\rm dof}$ denote the total number of multivariate basis functions $\Phi_{j,p}$. The spline space $\mathcal{V}_{h,p}$ can then be written as follows
\begin{eqnarray}
 \mathcal{V}_{h,p} = \text{span}\{ \Phi_{j,p} \mathbf{F}^{-1} \}_{j=1,\ldots,N_{\rm dof}}.
\end{eqnarray}

The Galerkin formulation of \eqref{eq:bilinearform} now becomes: Find $u_{h,p} \in \mathcal{V}_{h,p}$ such that 
\begin{equation} \label{eq:bilinearform2}
    a(u_{h,p},v_{h,p}) = (f_{h,p},v_{h,p}), \ \forall v_{h,p} \in V_{h,p}.
\end{equation}

The discretized problem in \eqref{eq:bilinearform2} can be written as a linear system
\begin{eqnarray} \label{eq:linearsystem1}
\left ( \mathbf{S}_{h,p} + k^2 \mathbf{M}_{h,p} - i k \mathbf{N}_{h,p} \right ) \mathbf{u}_{h,p} = \mathbf{f}_{h,p}.
\end{eqnarray}

Here, $\left (\mathbf{S}_{h,p}\right )_{i,j} = \int_{\Omega} \nabla \Phi_{i,p} \cdot \nabla \Phi_{j,p} \ \text{d}\Omega$ is the stiffness matrix, $\left (\mathbf{M}_{h,p} \right )_{i,j} = \int_{\Omega} \Phi_{i,p} \Phi_{j,p} \ \text{d}\Omega$ the mass matrix and $\left (\mathbf{N}_{h,p} \right)_{i,j} = \int_{\partial \Omega} \Phi_{i,p} \Phi_{j,p} \ \text{d}\Gamma$ the boundary mass matrix. Next, by defining $\mathbf{A}_{h,p} =\mathbf{S}_{h,p} + k^2 \mathbf{M}_{h,p} - i k \mathbf{N}_{h,p}$ we can write \begin{eqnarray} \label{eq:linearsystem}
 \mathbf{A}_{h,p} \mathbf{u}_{h,p} = \mathbf{f}_{h,p}.
\end{eqnarray}
For the ease of notation, we will proceed with the notation $\mathbf{Au} = \mathbf{f}$, and drop the subscript $(h,p)$. Using this discretization technique, we will now briefly explain the model problems used in this paper.

\section{Preconditioned Krylov Subspace Methods} \label{defgmres}
For Helmholtz-type problems, the number of degrees of freedom grows with the wave number $k$. Consequently, for larger values of $k$ the linear systems become very large, especially in two and three dimensions. As a result, direct solvers become unattractive and computationally expensive due to fill-in. Thus, in order to solve the model problems, an iterative method is considered. For normal matrices, the convergence of Krylov subspace methods is closely related to the underlying distribution of the eigenvalues. The more clustered the eigenvalues, the better and faster the method converges. 
For MP 1-B, we can easily deduce the analytical eigenvalues which are given by $\lambda_j = j^2\pi^2 - k^2$. It is easy to see that the resulting systems will have both positive and negative eigenvalues, rendering it indefinite. This limits our choice of Krylov subspace methods, where often GMRES is chosen as the underlying iterative solver. Many studies have investigated the performance of GMRES for the Helmholtz equation and the use of preconditioners is necessary in order to obtain satisfactory convergence. One of these preconditioners is the CSLP, which is defined by taking the original coefficient matrix $\mathbf{A}$ and adding a complex shift. Thus, in the one-dimensional case, CSLP $\mathbf{M}$ is given by 
\begin{equation} 
\mathbf{M} = \mathbf{A} +\beta_2 i k^2\mathbf{I},
\end{equation}
and the resulting preconditioned system becomes
\begin{equation} 
\mathbf{M^{-1}A} = \mathbf{M^{-1}f}.
\end{equation}
Here, $\mathbf{I}$ denotes the identity matrix and $\beta_2 \in \mathbb{R}$ the shift. In practice, the CSLP is often included by applying a fixed number of V-cycles of a (geometric) multigrid method to approximate $\textbf{M}^{-1}$. As a smoother within the multigrid method, we adopt damped Jacobi $(\omega = 0.6)$. Note that the use of standard smoothers (i.e. Jacobi or Gauss-Seidel) within a multigrid solver\cite{Gahalaut2013} in IgA results in $p-$dependent convergence. This has led to the development of non-standard smoothers to obtain $p$-independent convergence rates \cite{Hofreither2017a,Hofreither2017b,Donatelli2017,Sogn2019,Riva20,Tielen2020}. Their application within a multigrid method to approximate $\mathbf{M}^{-1}$ is, however, out of the scope of this paper. In order for $\mathbf{M}^{-1}$ to remain a good preconditioner, the shift $\beta_2$ should not be too small as otherwise multigrid will diverge \cite{gander2014optimized,ernst2012difficult}. On the other hand, the preconditioner should still remain close enough to the original coefficient matrix $\mathbf{A}$, which is also why $\beta_2$ should not be too large. \\
\\
While the complex shift transfers part of the unwanted spectrum onto the complex axis, unless the shift is kept very small, near-zero eigenvalues start appearing around the origin as the wave number increases \cite{erlangga2006,van2007spectral,sheikh2016accelerating}. This effect accumulates in higher-dimensions. Especially the real part of these near-zero eigenvalues is known to have a detrimental effect on the convergence behavior of the Krylov solver. One simple yet effective way to get rid of these unwanted near-zero eigenvalues is to use deflation. By using an orthogonal projection, the deflation operator, which we will denote by $\mathbf{P}$ projects these unwanted eigenvalues onto zero. Thus, for a general symmetric linear system, we can define the projection matrix $\mathbf{\widehat{P}}$ and its complementary projection $\mathbf{P}$ as
\begin{eqnarray}
 &\mathbf{\widehat{P}} = \mathbf{AQ} \mbox{ where } \mathbf{Q = ZE^{-1}Z^{T}} \mbox{ and } \mathbf{E} = \mathbf{Z^{T}AZ}, \label{adef1} \\ 
 &\mathbf{A} \, \in \mathbb{R}^{n \times n}, \, \mathbf{Z} \, \in \mathbb{R}^{m \times n}, \, \nonumber \\ 
 &\mathbf{P} = \mathbf{I - AQ}. \nonumber
\end{eqnarray}
Here the matrix $\mathbf{Z}$ is the deflation matrix whose columns consist of the deflation vectors and $\mathbf{E}$ denotes the coarse-grid variant of the original coefficent matrix $\mathbf{A}$. The performance of the deflation preconditioner depends on the choice of $\mathbf{Z}$. In principle, the deflation matrix is defined as the prolongation and restriction matrix from a multigrid setting using a first-order linear interpolation scheme \cite{sheikh110,sheikh112,gaul50,sheikh2014development,sheikh2016accelerating,garcia2018spectrum}. While this improves the convergence significantly, the near-zero eigenvalues start reappearing for very large wave numbers $k$. Consequently, it has been shown recently that the use of a quadratic interpolation scheme results in close to wave number independent convergence for the two-level deflation preconditioner \cite{dwarka2020scalable}. In fact, the use of these higher-order deflation vectors results in a smaller projection error compared to the case where a linear interpolation schemes is used. To construct the stencil for the deflation matrix $\mathbf{Z}$, we start by  introducing the rational $B\acute{e}zier$ curve. 
\begin{definition}[$B\acute{e}zier$ curve]
\label{bez1n}
A $B\acute{e}zier$ curve  of degree $n$ is a parametric curve defined by
\begin{eqnarray}
B(t) &= \sum\limits_{j=0}^{n} b_{j,n}(t) P_j, \hspace{2mm} 0 \leq t \leq 1, \hspace{2mm} \mbox{where the polynomials} \\
b_{j,n} (t) &= \left( n,j \right)  t^{j} (1 - t)^{n - j}, \hspace{2mm} j = 0,1, \hdots, n,
\end{eqnarray}
are known as the Bernstein basis polynomials of order $n$.
The points $P_j$ are called control points for the $B\acute{e}zier$ curve.
\end{definition}
\begin{definition}[Rational $B\acute{e}zier$ curve]
\label{bez2n}
A rational $B\acute{e}zier$ curve of degree $n$ with control points $P_0, P_1, \hdots, P_n$ and scalar weights $w_0, w_1, \hdots, w_n \in \mathbb{R}$ is defined as
\begin{equation}
C(t) = \frac{\sum\limits_{j=0}^{n} {w_j b_{j,n}(t)} {P_j}   }{\sum\limits_{j=0}^{n} {w_j b_{j,n}(t)}}.
\end{equation}
\end{definition}
For large $k$, the prolongation operator working on the even basis functions is not sufficiently accurate to map the underlying eigenvectors to its fine- and coarse-grid counterparts. We thus consider a quadratic rational $B\acute{e}zier$ curve in order to find appropriate coefficients to yield a higher order approximation of the fine-grid functions $u_{h}$ by the coarse grid functions $u_{2h}$. The motivation for using the rational $B\acute{e}zier$ curve is that the latter formulation allows for the weights to be adjusted in order to account for the higher requested accuracy at the even basis functions. In particular, if we define the coarse-grid basis function with respect to the degree of freedom $j$ by $[u_{2h}]_{j}$, then the quadratic approximation is defined as follows 
\begin{definition}[Quadratic Approximation]
\label{bez5n}
Let $[u_{2h}]_{({j-2})/2}$ and $[u_{2h}]_{({j+2})/2}$, be the neighbouring degrees of freedom of $[u_{2h}]_{j}$. Then the prolongation operator can be characterized by a Rational $B\acute{e}zier$ curve of degree 2 with polynomials
\begin{align*}
b_{0,2}(t) &= (1 - t)^2, \\
b_{1,2}(t) &= 2t(1 - t), \\
b_{2,2}(t) &= t^2,
\end{align*}
and $[u_{2h}]_{j/2}$, whenever $j$ is even. Because we wish to add more weight whenever $j$ is even, we take weights $w_0 = w_2 = \frac{1}{2}$, $w_1 = \frac{3}{2}$ and $t = \frac{1}{2}$ to obtain
\begin{align}
C(t) &= \frac{\frac{1}{2}{(1 - t)^2}{[u_{2h}]_{j-1}} + \frac{3}{2}{2t(1 - t)}{[u_{2h}]_{j}} + \frac{1}{2}{(t)^2}{[u_{2h}]_{j+1}}}{\frac{1}{2}{(1 - t)^2} + \frac{3}{2}{2t(1 - t)} + \frac{1}{2}{(t)^2}} \vspace{3mm} \nonumber \\
&= \frac{\frac{1}{2}(1 - \frac{1}{2})^2{[u_{2h}]_{j-1}} + \frac{3}{2}(2)(\frac{1}{2})(1 - \frac{1}{2}){[u_{2h}]_{j}} + \frac{1}{2}(\frac{1}{2})^2{[u_{2h}]_{j+1}}}{\frac{1}{2}(1 - \frac{1}{2})^2 + \frac{1}{2}(2)(\frac{1}{2})(1 - \frac{1}{2}) + \frac{1}{2}(\frac{1}{2})^2} \vspace{3mm} \nonumber \\
&= \frac{\frac{1}{8}{[u_{2h}]_{j-1}} + \frac{3}{4}{[u_{2h}]_{j}} + \frac{1}{8}{[u_{2h}]_{j+1}}}{1} \vspace{3mm} \nonumber \\
&=  \frac{1}{8}\left( [u_{2h}]_{j-1} + 6[u_{2h}]_{j} + [u_{2h}]_{j+1} \right). \nonumber \label{bez4}
\end{align}
When $j$ is odd, $[u_{2h}]_{({j-1})/2}$ and $[u_{2h}]_{({j+1})/2}$ are associated to an even degree of freedom and the resulting stencil leads to the standard linear interpolation scheme.
\end{definition}
Thus, with respect to the coarse-grid function $u_{2h}$ at degree of freedom $j$, we can define the stencil for the prolongation and restriction operator as
\begin{equation}
\left[\mathbf{Z} u_{2h}  \right ]_{j} = \left\{\begin{matrix}
\begin{matrix}
\frac{1}{8}  \left(\left[u_{2h}  \right ]_{\left(j - 2\right) / 2} + 6 \left[u_{2h}  \right ]_{\left(j \right) / 2} + \left[u_{2h}  \right ]_{\left(j + 2\right) / 2} \right) & \mbox{if $j$ is even,}\\
\frac{1}{2}  \left(\left[u_{2h}  \right ]_{\left(j - 1\right) / 2} + \left[u_{2h}  \right ]_{\left(j + 1\right) / 2} \right) &  \mbox{if $j$ is odd,}
\end{matrix}
\end{matrix}\right\},
\end{equation}
for $ j = 1, \dots, N_{\rm dof}$ and
\begin{align}
\left[\mathbf{Z^{T}} u_h  \right ]_{j} = \frac{1}{8}  \left(\left[u_h  \right ]_{\left( 2j - 2\right)} + 4  \left[u_h  \right ]_{\left( 2j + 1 \right)} + 6 \left[u_h  \right ]_{\left(2j \right)} + 4 \left[u_h  \right ]_{\left( 2j + 1 \right)} + \left[u_h  \right ]_{\left( 2j + 2\right)} \right),
\end{align}
for $j = 1, \dots, \frac{N_{\rm dof}}{2}$. 
Now that we have a stencil to construct $\mathbf{Z}$, we can use Equation \eqref{adef1} to construct the deflation preconditioner. The resulting linear system to be solved becomes 
\begin{eqnarray}
\mathbf{P^{T}A} = \mathbf{P^{T}f}.
\end{eqnarray}
Often, the deflation preconditioner $\mathbf{P}$ is combined with the CSLP $\mathbf{M}$ to accelerate convergence, which leads to solving the following system
\begin{eqnarray}
\mathbf{P^{T}M^{-1}Au} &= \mathbf{(I - AQ)^{T}M^{-1}Au}
&= \mathbf{(I - AQ)^{T}M^{-1}f},
\end{eqnarray}
where, as mentioned previously, $\mathbf{M^{-1}}$ is generally approximated using a multigrid method. Note that the operator $\mathbf{P^{T}}$ is never constructed explicitly but is passed as a function handle onto the coefficient matrix $\mathbf{A}$ within the GMRES-algorithm. Moreover, we will refer to $\mathbf{P}$ based on the higher-order quadratic approximation as the 'Adapted Deflation Preconditioner' (ADP) to distinguish between the standard deflation preconditioner using linear interpolation and the higher-order deflation scheme. Additionally, a weight-parameter can be included to further increase the accuracy of the prolongation and restriction operator \cite{dwarka2020scalable}. In this case, the stencil for the prolongation and restriction operator is given by 
\begin{equation}
\left[\mathbf{Z} u_{2h}  \right ]_{j} = \left\{\begin{matrix}
\begin{matrix}
\frac{1}{8}  \left(\left[u_{2h}  \right ]_{\left(j - 2\right) / 2} + (6 - \varepsilon) \left[u_{2h}  \right ]_{\left(j \right) / 2} + \left[u_{2h}  \right ]_{\left(j + 2\right) / 2} \right) & \mbox{if $j$ is even,}\\
\frac{1}{2}  \left(\left[u_{2h}  \right ]_{\left(j - 1\right) / 2} + \left[u_{2h}  \right ]_{\left(j + 1\right) / 2} \right) &  \mbox{if $j$ is odd,}
\end{matrix}
\end{matrix}\right\},
\end{equation}
for $ j = 1, \dots, N_{\rm dof}$ and
\begin{align}
\left[\mathbf{Z^{T}} u_h  \right ]_{j} = \frac{1}{8}  \left(\left[u_h  \right ]_{\left( 2j - 2\right)} + 4  \left[u_h  \right ]_{\left( 2j + 1 \right)} + (6 - \varepsilon) \left[u_h  \right ]_{\left(2j \right)} + 4 \left[u_h  \right ]_{\left( 2j + 1 \right)} + \left[u_h  \right ]_{\left( 2j + 2\right)} \right),
\end{align}
for $j = 1, \dots, \frac{N_{\rm dof}}{2}$. Note that the value of $\varepsilon$ is constant with respect to $k$ and $kh$ and is chosen such that the projection error is minimized \cite{dwarka2020scalable}.

% \begin{enumerate}
%     \item Brief introduction to Krylov Subspace methods with focus on GMRES. - check
%     \item Preconditioning techniques for the Helmholtz equation which will be used in this paper. - check
%     \begin{enumerate}
%         \item CSLP - check
%         \item Deflation - check
%     \end{enumerate}
%     \item Introduction and motivation for using deflation methods. - check
%     \item Explain construction of coarse grid space/linear system. 
%     \item Give scheme for the Quadratic Bezier Curve. 
%     \item Give an overview of the various linear systems to be solved.
% \end{enumerate}

\section{Numerical Results} \label{numres}
To assess the quality of the proposed iterative solver, we consider the model problems described in \ref{variational}. We start by studying the pollution error for our one-dimensional model problem when adopting high-order B-spline basis functions for the spatial discretization. In \cite{diwan2020iterative}, a detailed first application of IgA discretizations for Helmholtz problems has been given. We therefore only show the pollution reduction for the model problems used in this paper. 
We proceed by conducting a spectral analysis in one dimension (MP 1-B) to investigate the effect of the proposed preconditioning techniques on the spectrum of the preconditioned operator. Finally, the convergence of the iterative solver is studied in terms of 
both iteration numbers and CPU timings. These are obtained for the proposed deflation based preconditioner and compared to the use of the (exactly inverted) CSLP. 

\subsection{Pollution Error}
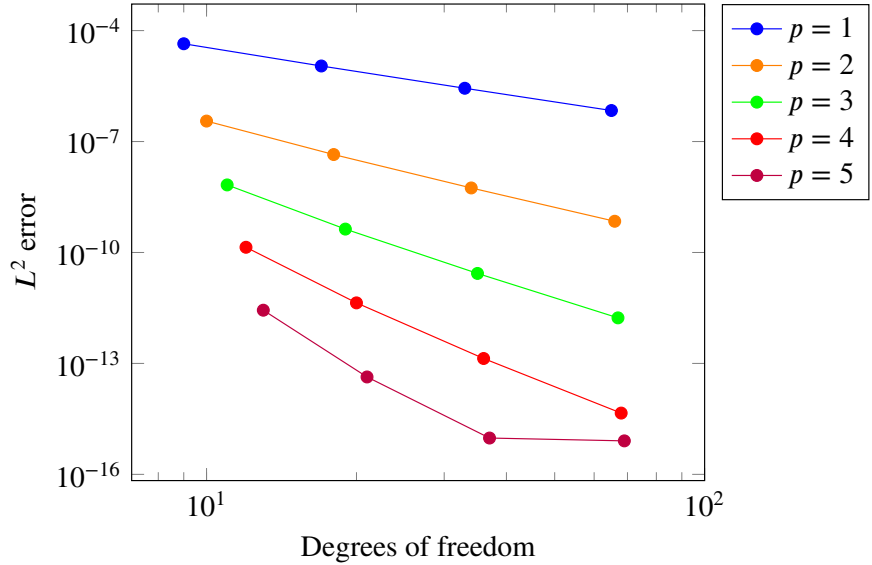
\begin{figure}[ht!]
\centering
   \begin{tikzpicture}[scale=1.1]
    \begin{axis}[xmode=log,xlabel=Degrees of freedom, xmax=100,ymode=log,ylabel=$L^2$ error, legend pos = outer north east]
       \addplot[blue,mark=*] plot coordinates 
      {
        (9,    4.409E-5)
        (17,   1.104E-5)
        (33,   2.760E-6)
        (65,   6.900E-7)
        %(129,  1.725E-7)
        %(257,  4.313E-8)
        %(513,  1.078E-8)
        };
      \addlegendentry{$p=1$}
        \addplot[orange,mark=*] plot coordinates 
      {
        (10,    3.576E-7)
        (18,   4.452E-8)
        (34,   5.560E-9)
        (66,   6.986E-10)
        %(130,  8.684E-11)
        %(258,  1.086E-11)
        %(514,  1.357E-12)
        };
      \addlegendentry{$p=2$}
       \addplot[green,mark=*] plot coordinates 
      {
        (11,    6.712E-9)
        (19,   4.289E-10)
        (35,   2.713E-11)
        (67,   1.706E-12)
        %(131,  1.069E-13)
        %(259,  2.151E-14)
        %(515,  1.003E-14)
        };
      \addlegendentry{$p=3$}
      \addplot[red,mark=*] plot coordinates 
      {
        (12,    1.383E-10)
        (20,   4.336E-12)
        (36,   1.359E-13)
        (68,   4.536E-15)
        %(132,  2.328E-14)
        %(260,  1.182E-13)
        %(516,  5.242E-13)
        };
      \addlegendentry{$p=4$}
      \addplot[purple,mark=*] plot coordinates 
      {
        (13,    2.747E-12)
        (21,   4.303E-14)
        (37,   9.594E-16)
        (69,   8.045E-16)
        %(133,  1.354E-14)
        %(261,  1.661E-14)
        %(517,  1.557E-13)
        };
      \addlegendentry{$p=5$}
      \end{axis}
\end{tikzpicture}
\caption{Spatial convergence for different values of $p$ obtained with (deflated) GMRES for MP 1-A, where $k=1$.}
\label{fig:1}
\end{figure}

As a first verification of the quality of the solver, a spatial convergence test has been performed for the MP 1-A benchmark for a fixed value of the wave number ($k=1$). Figure \ref{fig:1} shows the $L_2$-error under mesh refinement for different values of $p$ obtained with a (deflated) GMRES solver. Note that, for all values of $p$, the order of convergence observed is $\mathcal{O}(h^{p+1})$, as expected from literature \cite{hughes2005}. For $p=5$ and a sufficiently fine mesh, the $L_2$-error becomes close to machine precision and therefore suffers from errors in floating point operations. Detailed $L_2$-errors can be found in Table \ref{tab:1}.

\begin{table}[ht!]
\centering
\scalebox{1.00}{
\begin{tabular}{c|c|c|c|c|c|}
    $N_{\rm dof}-p$  &     $p=1$               &     $p=2$                &     $p=3$                &     $p=4$                 &     $p=5$ \\ \hline
    $8$       &    $4.409\cdot10^{-5}$  &     $3.576\cdot10^{-7}$  &     $6.711\cdot10^{-9}$  &     $1.383\cdot10^{-10}$  &     $2.747\cdot10^{-12}$ \\
    $16$      &    $1.104\cdot10^{-5}$  &     $4.452\cdot10^{-8}$  &     $4.289\cdot10^{-10}$ &     $4.336\cdot10^{-12}$  &     $4.303\cdot10^{-14}$ \\
    $32$      &    $2.760\cdot10^{-6}$  &     $5.560\cdot10^{-9}$  &     $2.713\cdot10^{-11}$ &     $1.359\cdot10^{-13}$  &     $9.594\cdot10^{-16}$ \\
    $64$      &    $6.900\cdot10^{-7}$  &     $6.948\cdot10^{-10}$ &     $1.706\cdot10^{-12}$ &     $4.536\cdot10^{-15}$  &     $8.045\cdot10^{-16}$ \\
\end{tabular}
}
\caption{$L_2$-error under mesh refinement for different values of $p$ obtained with (deflated) GMRES for MP 1-A, where $k=1$.}
\label{tab:1}
\end{table}
\FloatBarrier
In order to determine the effect of using B-spline basis functions on the pollution error, we present the $L_2$-error as a function of the wave number $k$ as well. Note that the case $p=1$ corresponds to the standard Lagrangian FEM solution. We observe that for $p=2$ to $p=5$ the $L_2$-error with respect to the analytical solution decreases. While this leads to significant more accurate solutions, we do observe that as the wave number increases, the $L_2$-error increases accordingly. This is in line with the literature, as it has been proven that the pollution error can not be avoided completely \cite{pollution1997,pollution1999}. Moreover, as $k$ increases the advantage of using $p = 5$ over $p = 4$ decreases as both lead to similar accuracy. For standard FEM, this was already observed \cite{singer1998high}. Furthermore, decreasing the number of degrees of freedom per wavelength from 10 (solid line) to 7.5 (dashed line) already results in lower accuracy. In fact, the achieved accuracy for $p=4$ and $p=5$ with 7.5 degrees of freedom per wavelength is similar to the obtained accuracy for $p = 3$ when 10 degrees of freedom per wavelength are used. Thus, in order to warrant for sufficiently accurate numerical solutions for larger wave numbers, we will keep the grid resolution at $kh = 0.625$.

\begin{figure}[ht!]
\centering
   \begin{tikzpicture}[scale=1.1]
    \begin{axis}[xmode=log,xlabel=Wave number, ymode=log,ylabel=$L^2$ error,legend pos = outer north east, legend columns = 2]
       \addplot[blue,mark=*] plot coordinates 
      {
        (100,    0.026506685961951)
        (500,   0.0413343411140983)
        (1000,   0.0440869399374673)
        (2500,   0.0437081856179975)
        (5000,  0.044095191495)
        (7500,  0.0443810551015249)
        (10000,  0.044401180836523)
        };
      \addlegendentry{$p=1, kh =0.625$}
      
      \addplot [blue,dashed,mark=*] plot coordinates
      {
        (101,    0.0401296837517977)
        (501,    0.04276786568326183)
        (1001,   0.0430883650044422)
        (2501,   0.0440919641365801)
        (5001,  0.0437701233674736)
        (7501,  0.0440215609146077)
        (10001,   0.0438586842256983)
        };
      \addlegendentry{$p=1, kh =0.825$}
      \addplot[orange,mark=*] plot coordinates 
      {
        (102,   0.000221419686435354)
        (502,   0.00105570626590273)
        (1002,   0.00210417362251059)
        (2502,   0.00524835773967353)
        (5002,  0.0104160905411237)
        (7502,  0.0154586281683434)
        (10002,  0.0203127746872357)
        };
      \addlegendentry{$p=2, kh =0.625$}
       \addplot [orange,dashed,mark=*] plot coordinates
      {
        (103,    0.000740723681328122)
        (503,    0.00349052229478681)
        (1003,   0.0069346442550872)
        (2503,   0.0169298040698228)
        (5003,   0.031621400925704)
        (7503,   0.0423474814608057)
        (10003,  0.0481742133296658)
        };
        \addlegendentry{$p=2, kh =0.825$}
       \addplot[green,mark=*] plot coordinates 
      {
        (104,   6.20008814284829e-06)
        (504,   1.47991771734765e-05)
        (1004,   2.79612273185094e-05)
        (2504,   6.81957909553223e-05)
        (5004,  0.000136497576683339)
        (7504,  0.000204670024111978)
        (10004,  0.000272830729907881)
        };
      \addlegendentry{$p=3, kh =0.625$}      
      \addplot [green,dashed,mark=*] plot coordinates
      {
        (101,    2.74345900207248e-05)
        (501,    8.75217829852007e-05)
        (1001,   0.000170342382306206)
        (2501,   0.000419157800004363)
        (5001,   0.000836331049466228)
        (7501,   0.00125348174119513)
        (10001,  0.00167054725576272)
        };
       \addlegendentry{$p=3, kh =0.825$}
       \addplot[red,mark=*] plot coordinates 
      {
        (103,   1.12405076192994e-06)
        (503,   4.28533316254618e-06)
        (1003,   8.41519378475062e-06)
        (2503,   2.06537135289553e-05)
        (5003,  4.15644686116019e-05)
        (7503,  6.23339981438306e-05)
        (10003,  8.30932245116405e-05)
        };
      \addlegendentry{$p=4, kh =0.625$}      
      \addplot [red,dashed,mark=*] plot coordinates
      {
        (104,    6.35322905104563e-06)
        (504,    2.39204864072859e-05)
        (1004,   4.68763050821601e-05)
        (2504,   0.000115432998046555)
        (5004,   0.00023023644416439)
        (7504,   0.00034502075057239)
        (10004,  0.000459799682009423)
        };
      \addlegendentry{$p=4, kh =0.825$}
      \addplot[purple,mark=*] plot coordinates 
      {
        (104,  1.20272093008407e-06)
        (504,  1.93916104585525e-06)
        (1004,  3.70212011576164e-06)
        (2504,  9.1871606775995e-06)
        (5004, 1.83482939392166e-05)
        (7504, 2.7511946848068e-05)
        (10004, 3.66747834845872e-05)
        };
      \addlegendentry{$p=5, kh =0.625$}      
      \addplot [purple,dashed,mark=*] plot coordinates
      {
        (104,    4.40948759309054e-06)
        (504,     1.15360587640238e-05)
        (1004,   2.24991208299944e-05)
        (2504,   5.54944770252471e-05)
        (5004,   0.000110650327300312)
        (7504,   0.000165810686078885)
        (10004,  0.000220971362445622)
        };
      \addlegendentry{$p=5, kh =0.825$}
      \end{axis}
\end{tikzpicture}
    \caption{L2-error for MP 1-A using $p=1$ to $p=5$ for various wave numbers $k$. The grid resolution has been set to $kh \approx 0.625$ (solid) and $kh = 0.825$ (dashed).}
\label{fig:1v}
\end{figure}
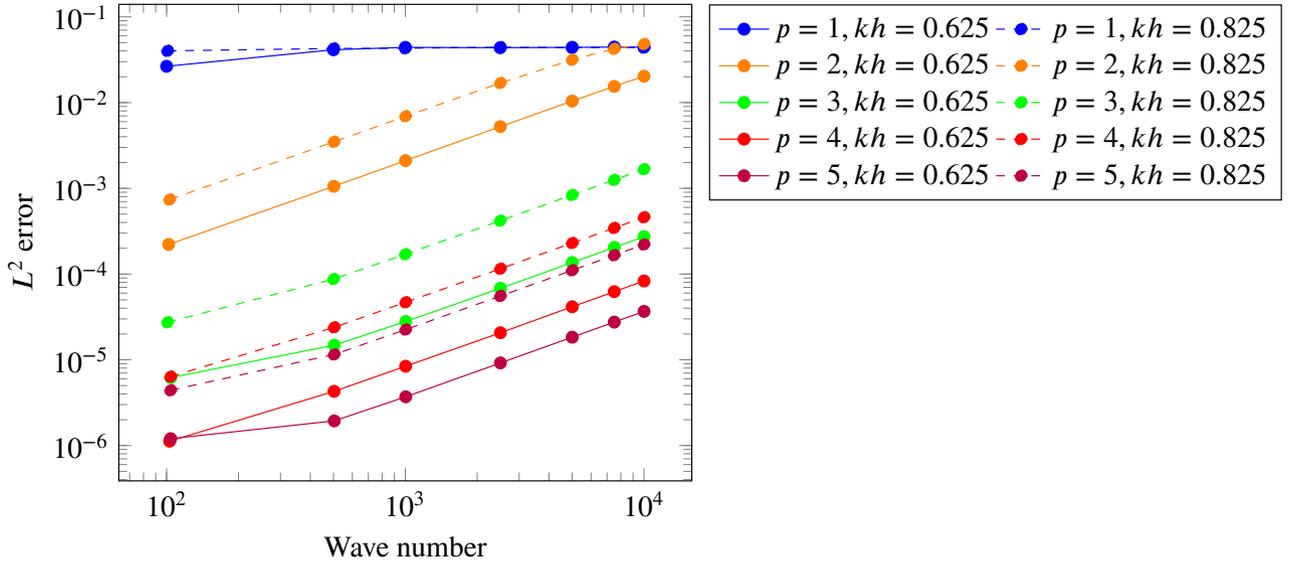
\FloatBarrier

\subsection{Spectral Analysis} \label{spec}
\FloatBarrier
We now proceed by analyzing the spectrum of the preconditioned system of MP 1-B. 
It is widely known that the near-zero eigenvalue distribution strongly affects the resulting convergence factor of Krylov subspace methods. In general, these eigenvalues close to the origin hamper the convergence of such methods. By using Dirichlet boundary conditions, we additionally have the most unfavorable distribution of eigenvalues, allowing us to fully examine the potency of the preconditioner. 
With respect to CSLP, many studies have confirmed that unless the complex shift is kept very small and the inversion is performed exactly, the eigenvalues cluster near the origin \cite{van2007spectral,gander2014optimized,gander2015applying}. 
In this work, we are not inverting the CSLP exactly and we thus need to derive a proxy of the multigrid iteration used to approximate the inverse. This can be done by using the two-grid iteration matrix from a multigrid setting \cite{erlangga2006comparison}. This leads to the following approximation for $\mathbf{M}$
\begin{eqnarray*}
\tilde{\mathbf{M}}^{-1} \approx \left ( \mathbf{I} - (\omega\mathbf{D})^{-1} \mathbf{M} \right )^{\nu} \left ( \mathbf{I} - \mathbf{Z} \mathbf{M}_{2h}^{-1} \mathbf{Z}^{\top} \right ) \left ( \mathbf{I} - (\omega\mathbf{D})^{-1} \mathbf{M} \right )^{\nu}, 
\end{eqnarray*}
where $\mathbf{M}_{2h}$ denotes the coarse-grid variant of the CSLP, $\mathbf{D}$ the diagonal of $\mathbf{M}$ and $\nu$ denotes the smoothing steps. Additionally, we use damped Jacobi as a smoother with damping parameter $\omega=0.6$. Note that for the multigrid cycle, $\mathbf{Z}$ is now the standard geometric multigrid prolongation and restriction operator based on the linear interpolation scheme. 
Using this approximation for $\textbf{M}^{-1}$, we study the eigenvalues of the linear system $\mathbf{P^{T}}\mathbf{\tilde{M}}^{-1}\mathbf{A}$, where $\mathbf{P}$ denotes the adapted deflation preconditioner based on the quadratic Bezier scheme. 

\begin{figure}[ht!]
    \centering
    \includegraphics[scale=0.41]{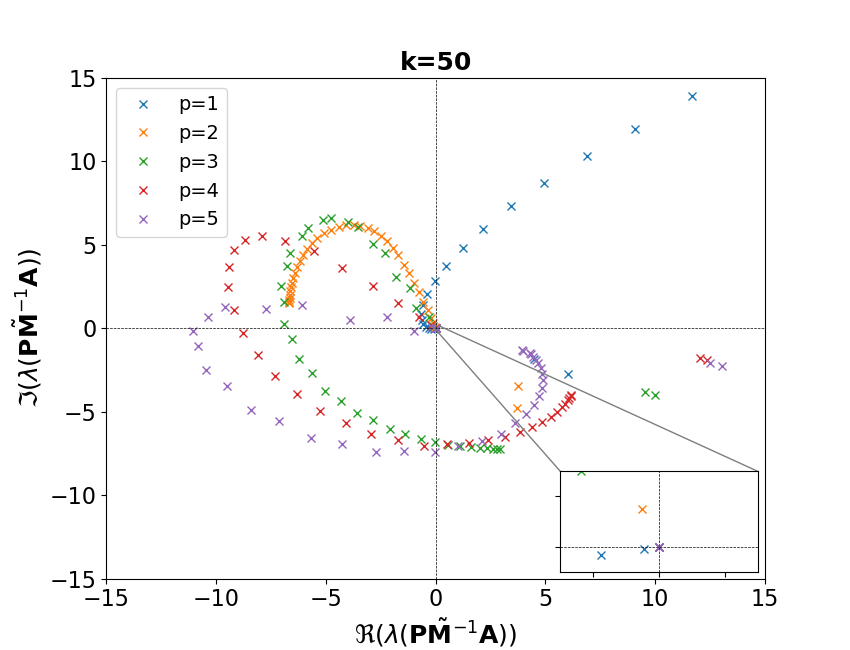}
    \includegraphics[scale=0.41]{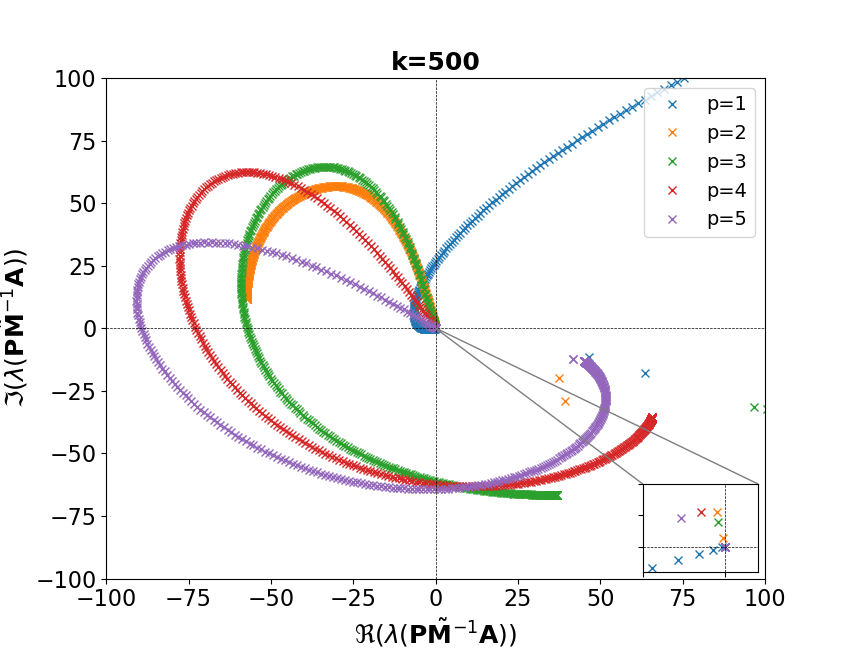}
    \caption{Spectrum of the preconditioned operator $\mathbf{P}\tilde{\mathbf{M}}^{-1}\mathbf{A}$ for different values of $p$, where $k=50$ (left) and $k=500$ (right) for MP 1-A.}
    \label{fig:K}
\end{figure}

Figure \ref{fig:K} shows the spectra of the preconditioned linear system for $k=50$ (left) and $k=500$ (right) for different values of $p$. The complex shift has been set to $\beta_2=1$ and one pre- and post-smoothing step has been used. Note that half of the eigenvalues of the preconditioned system will be projected onto the origin. The other half of the eigenvalues will therefore be non-zero. 
% {\color{blue} Looking at the non-zero eigenvalues for $k = 50$ (left), we can see that for all values of $p$, most eigenvalues are located away from the origin. However, as $p$ increases the number of near-zero eigenvalues with respect to the real-axis starts to grow. Once we increase the wave number to $k=500$, we detect some more near-zero eigenvalues, especially for $p=1$ (blue). Therefore, we expect the number of GMRES iterations to have a slight dependence on the approximation order $p$ of the B-spline basis functions. 
% If we focus on the small box in the figures containing a detailed illustration of what is occurring near the origin, we observe that the visual angle these near-zero eigenvalues make with respect to the real-axis appears to be fairly constant for both $k = 50$ and $k=500$. Thus, while we do expect slightly more iterations when $p$ increases, this will be at a lower rate compared to the generally observed $k-$dependent convergence for multigrid approximated CSLP without deflation. In the latter case, the number of iterations grows fast as the wave number $k$ increases.} \\
% {\color{cyan} 
For $k=50$ (left), all eigenvalues for a fixed value of $p$ have a spiral shape, apart for the case $p=1$. Furthermore, the angle between the eigenvalues and the real-axis in Quadrant $2$ becomes smaller for higher values of $p$. Therefore, we can expect a $p$-dependency for small values of $k$ for $p \geq 2$. For $k=500$ this becomes even more obvious visually, as the higher number of degrees of freedom leads to more eigenvalues. As the preconditioned operator becomes too large to determine all eigenvalues, it remains unsure how the spectra will further developed for large values of $k$.

% As a result, hat, when applying the CSLP exactly, the preconditioned linear system can written as $\mathbf{P}\mathbf{M}^{-1}\mathbf{A}$. Here $\mathbf{P}$ denotes the deflation operator and $\mathbf{M}^{-1}$ the CSLP preconditioner. However, as we do not invert the preconditioner exactly, but apply a multigrid method to approximate the inverse, $\mathbf{M}^{-1}$ is replaced by $\tilde{\mathbf{M}}^{-1}$. Adopting (damped) Jacobi as a smoother within a two-level multigrid method, $\tilde{\mathbf{M}}^{-1}$ can be written as:

% \begin{eqnarray}
% \tilde{\mathbf{M}}^{-1} = \left ( \mathbf{I} - (\omega\mathbf{D})^{-1} \mathbf{M} \right )^{\nu} \left ( \mathbf{I} - \mathbf{Z} \mathbf{A} \mathbf{Z}^{\top} \right ) \left ( \mathbf{I} - (\omega\mathbf{D})^{-1} \mathbf{M} \right )^{\nu}. \nonumber 
% \end{eqnarray}

% Figure \ref{fig:K} shows the spectra of the preconditioned linear system for $K=50$ (left) and $K=500$ (right) for different values of $p$. For all values of $p$, most eigenvalues are located away from the origin, but some remain close to zero, in particular for $p=5$. Therefore, it is expected that the number of GMRES iterations when solving the resulting linear system could slightly depend on the approximation order $p$ of the B-spline basis functions.  

% {\color{blue} Next, in Figure \ref{fig:p}, we fix $p=2$ (left) and $p=5$ (right) and let $k$ increase from $k=50$ to $k = 250$. Here we can clearly observe that for $p=2$ eigenvalues remain fairly clustered. Increasing $k$ leads to a larger spread between the eigenvalues. Increasing $k$ primarily leads to the imaginary part increases as we see the eigenvalues moving upwards. If we focus on the small box containing a detailed illustration of what is occurring near the origin, we observe that for larger $k$ more and more eigenvalues are starting to move closer towards the origin. Closest to the origin we can clearly see the eigenvalues for $k  = 250$ (purple) and $k = 200$ (red) appearing. 
% If we now take $p=5$, we observe a similar effect, however the eigenvalues seem less clustered when compared to the case $p =2$. An interesting observation can be made with respect to $p=5$. As $k$ increases, the imaginary part of the eigenvalues starts increasing and the eigenvalues in Quadrant II are moving away from the real-axis. If we focus on the detailed view near the origin, we observe smaller eigenvalues for larger $k$ which was also reported for the case $p = 2$. However, comparing $p = 2$ to $p = 5$, note that as $k$ increases, more near-zero eigenvalues are clustering near the origin when $p = 2$.  Thus, the larger the wave number, the more near-zero eigenvalues we expect for smaller $p$. We thus expect to see an inverse relation between the number of iterations and larger $p$ for very large wave numbers. }

\begin{figure}[ht!]
    \centering
    \includegraphics[scale=0.41]{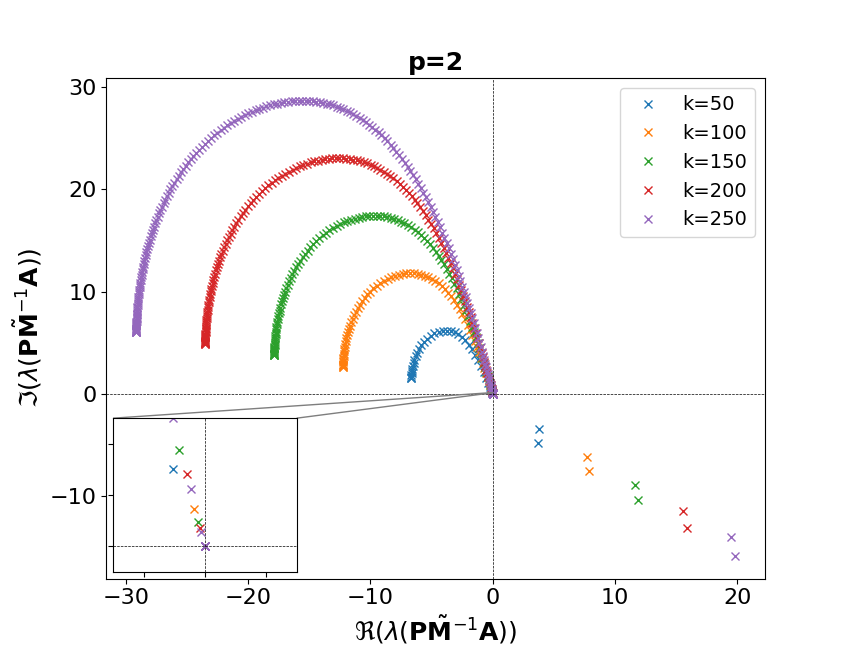}
    \includegraphics[scale=0.41]{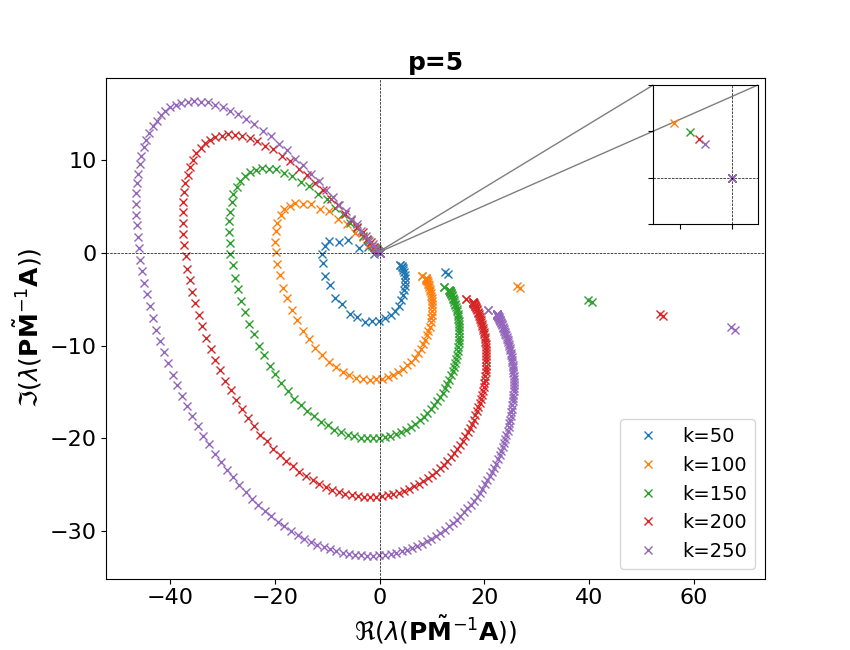}
    \caption{Spectrum of the preconditioned operator $\mathbf{P}\tilde{\mathbf{M}}^{-1}\mathbf{A}$ for different values of $k$, where $p=2$ (left) and $p=5$ (right) for MP 1-A. No weight-parameter has been included.}
    \label{fig:p}
\end{figure}

Next, in Figure \ref{fig:p}, we fix $p=2$ (left) and $p=5$ (right) and let $k$ increase from $k=50$ to $k = 250$. Here we can clearly observe that for $p=2$, the eigenvalues remain fairly clustered in a semi-circular shape. Increasing $k$ leads to a larger radius of this semi-circle and therefore a larger spread of the eigenvalues. If we focus on the small box containing a detailed illustration of what is occurring near the origin, we observe that for larger $k$ more and more eigenvalues are starting to move closer towards the origin. Closest to the origin we can clearly see the eigenvalues for $k  = 250$ (purple) and $k = 200$ (red) appearing. Although the eigenvalues seem less clustered for $p=5$, the same general behavior can be observed.

\begin{figure}[ht!]
    \centering
    \includegraphics[scale=0.41]{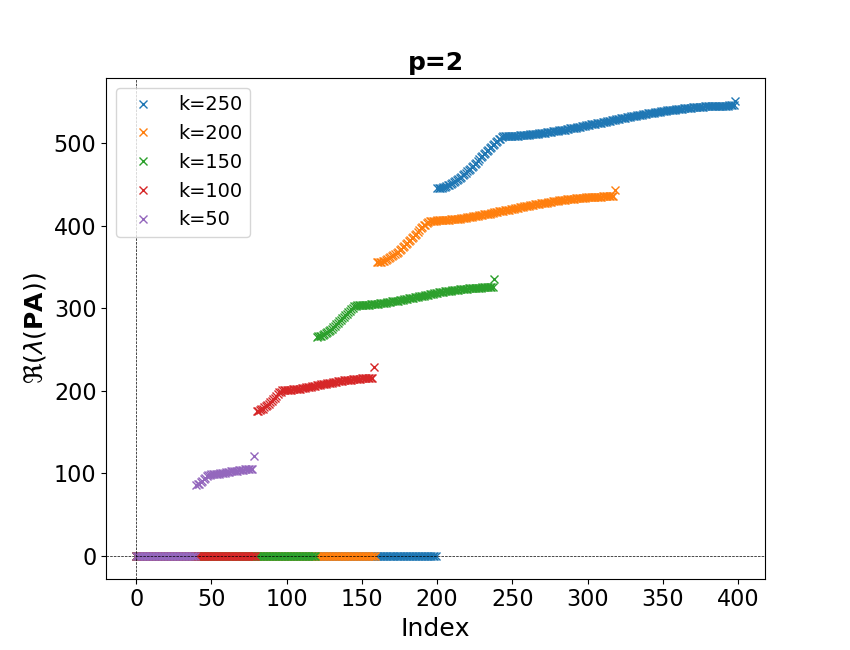}
    \includegraphics[scale=0.41]{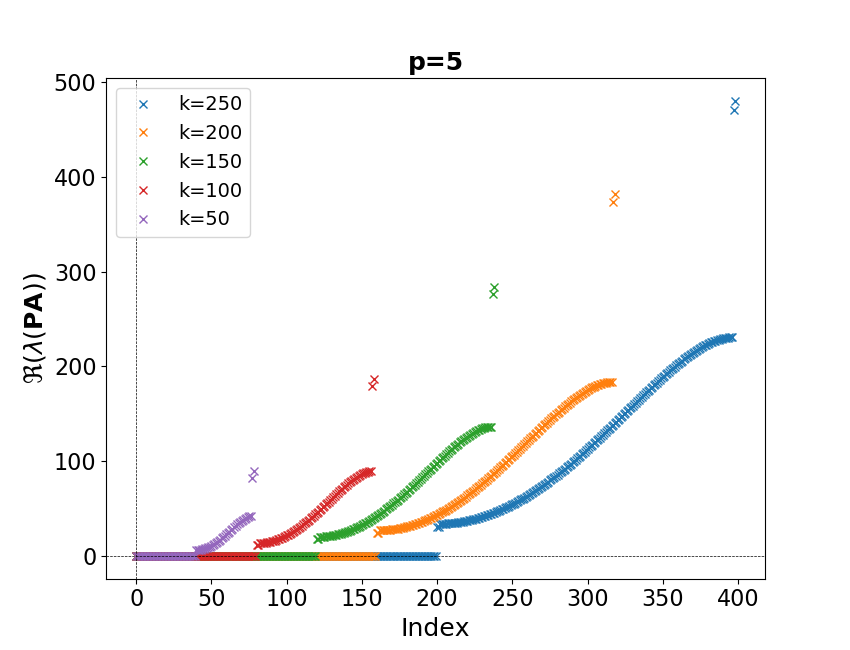}
    \caption{Real part of the spectrum of the preconditioned operator $\mathbf{P}\mathbf{A}$ for different values of $k$, where $p=2$ (left) and $p=5$ (right) for MP 1-A. No weight-parameter has been included.}
    \label{fig:pa1}
\end{figure}

Classically, deflation based preconditioners are combined with the CSLP in order to obtain faster GMRES-convergence. Note that the projection matrix $\mathbf{P}$ projects a certain part of the spectrum of the coefficient matrix $\mathbf{A}$ onto zero. The addition of the CSLP ensures that the remaining non-zero eigenvalues are shifted towards the complex axis, which gives it the typical circular spectrum in the complex plane. However, for finite differences discretizations, the use of the CSLP is often redundant as wave number independent convergence can already be attained by using deflation without another preconditioner. An interesting point of investigation would be to study the spectrum of the preconditioned system $\mathbf{PA}$. 
In Figure \ref{fig:pa1}, we study the spectrum of $\mathbf{PA}$ where we use the weight-parameter $\varepsilon$ in order to construct accurate higher-order deflation vectors. 
We indeed observe that half of the eigenvalues are mapped onto zero and the remaining part of the eigenvalues remains clustered. The eigenvalues no longer cross the negative real axis, which results in the preconditioned system $\mathbf{PA}$ being positive semi-definite. Apart from a scaling factor, the spectrum of $k = 50$ looks similar to the spectrum of $k = 250$ and illustrative of the $k-$independent convergence. However if we compare $p = 2$ (left) to $p = 5$ (right), we observe that for $p = 5$ the eigenvalues of $\mathbf{PA}$ are closer to zero and have a larger spread between the smallest and largest eigenvalue. For example for $k = 250$, the eigenvalues for $p = 2$ lie in the ballpark of 450 to 550, whereas for $p = 5$ the eigenvalues lie between 50 and 250.

\begin{figure}[ht!]
    \centering
    \includegraphics[scale=0.41]{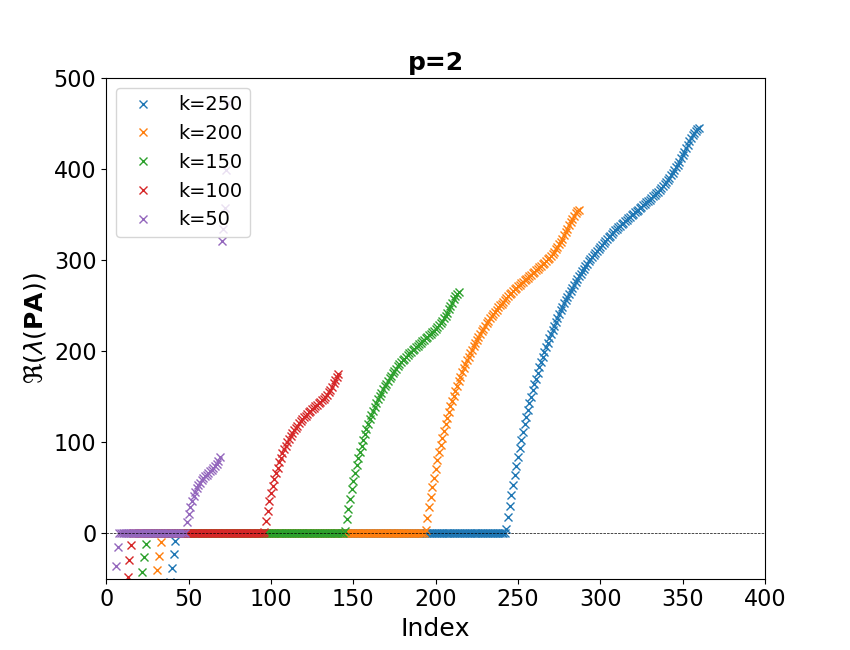}
    \includegraphics[scale=0.41]{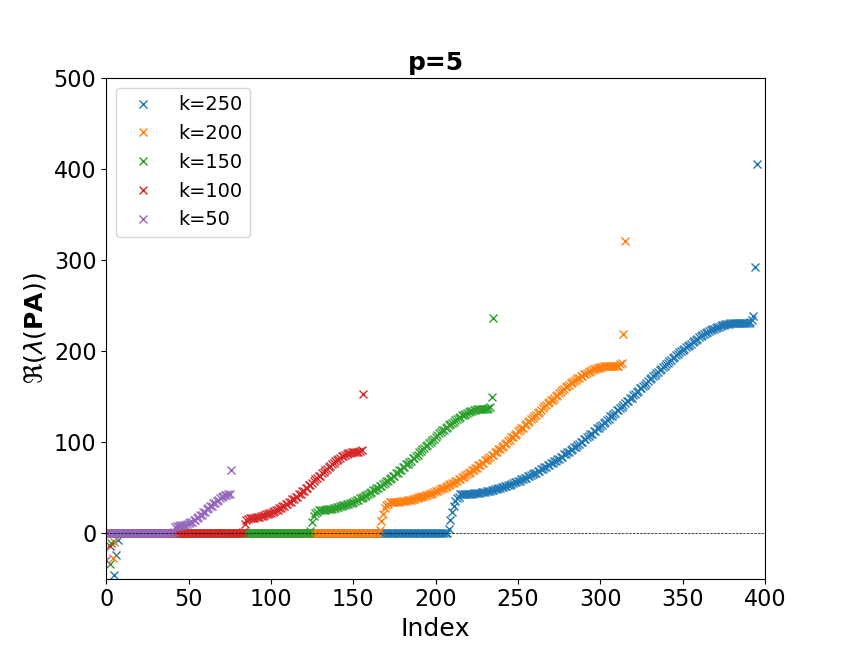}
    \caption{Real part of the spectrum of the preconditioned operator $\mathbf{P}\mathbf{A}$ for different values of $k$, where $p=2$ (left) and $p=5$ (right) for MP 1-A. Here we have used the weight-parameter $\varepsilon$.}
    \label{fig:pa2}
\end{figure}

For illustration purposes, we study the effect of interpolating and restricting the fine-grid systems with low accuracy. In Figure \ref{fig:pa2}, we have plotted the spectrum of $\mathbf{PA}$, where we deliberately set the weight-parameter to a value which lowers the accuracy of the interpolation scheme to construct the deflation matrix $\mathbf{Z}$. It immediately becomes apparent that the resulting preconditioned system is again indefinite as some eigenvalues are still negative. Moreover, if we compare $p =2$ (left) to $p =5$ (right), we observe a larger spread for $p=2$ compared to $p=5$. This is the opposite of what we observed in Figure \ref{fig:pa1}. In both cases, the example is illustrative of the fact that having a low-order interpolation scheme to construct the prolongation and restriction operator, will lead to an ineffective mapping of the underlying eigenvalues and eigenvectors. As the wave number increases and the solutions become more oscillatory, the accurate mapping of the fine- and coarse-space become of increasing importance. Therefore, we chose a weight-parameter such that the projection error with respect to the eigenvectors are minimized \cite{dwarka2020scalable}.
\FloatBarrier

\subsection{Numerical experiments}
We will now present the convergence results for our model problems using the preconditioners described above. Unless stated otherwise, we set the grid resolution at $kh \approx 0.625$, which is equivalent to using 10 degrees of freedom per wavelength. We use GMRES as the underlying Krylov subspace method and use a stopping criterium on the relative residual of $10^{-7}$. A serial implementation is considered on an Intel(R) i7-8665 CPU @ 1.90GHz using 8GB of RAM. \\
For the sake of completeness and clarity, we briefly introduce the notation of the preconditioners used in the experiments. 
\begin{itemize}
    \item $D$ := Adapted Deflation Preconditioner (ADP) + GMRES. 
     \item $D_{\varepsilon}$ := Adapted Deflation Preconditioner (ADP) + GMRES using the shift-parameter $\varepsilon$ to construct the deflation matrix. The value has been taken from \cite{dwarka2020scalable} and is constant throughout the use of the numerical experiments. 
    \item $C_{ex}$ := CSLP (exactly inverted) + GMRES.
    \item $DC_{MG}^{j}$ := $ADP$ preconditioner + GMRES using $j$ number of multigrid V-cycles combined with (damped) Jacobi smoothing. 
    \item $D_{\varepsilon}C_{MG}^{j}$ := $AD_{\varepsilon}P$ preconditioner + GMRES using $j$ number of multigrid V-cycles combined with (damped) Jacobi smoothing. 
\end{itemize}

\subsubsection{One-dimensional model problems}
\paragraph{MP 1-B}
\FloatBarrier
We start by numerically solving MP 1-B using the deflation preconditioner together with the multigrid approximation of the CSLP. We differentiate between deflation with and without the weight-parameter $\varepsilon$ and we vary the number of V-cycles between 1 and 10 iterations to obtain a fair approximation of the inverse of the CSLP. Table \ref{tab:4} shows the number of GMRES iterations for the three different combinations. Starting with $DC_{MG}^{1}$ (first column) we observe that the number of iterations both grow with $k$ and $p$. These results are in line with the spectral analysis from Section \ref{spec}, in particular Figure \ref{fig:K} and Figure \ref{fig:p}. There we observed that the angle the eigenvalues make with the real axis becomes smaller for increasing $p$, anticipating some $p-$dependent convergence. Similarly, in Figure \ref{fig:K}, the radius of the circular shape of the eigenvalues grows with $k$, leading to the expectation that the number of iterations could grow with $k$. However, for very large wave numbers such as $k = 10^4$, we observe that the number of iterations is inversely related to $p$. Note that the spectrum of such large wave numbers has not been examined in this work. \\
For $D_{\varepsilon}C_{MG}^{1}$ (second column) we solely observe $p-$dependent convergence. Once we add the weight parameter $\varepsilon$ to the deflation preconditioner we obtain $k$-independent convergence up to $10^6$. Finally, increasing the number of V-cycles to 10 for $D_{\varepsilon}C_{MG}^{10}$ (third column) leads to $p$-independent convergence and shows identical results to inverting CSLP exactly; see Table \ref{tab:2} and Table \ref{tab:3}. Note, however, that application of $D_{\varepsilon}C_{MG}^{10}$ is more expensive compared to the application of $D_{\varepsilon}C_{MG}^{1}$ as we use more V-cycles in order to obtain a fair approximation of the CSLP. This result, however, is in line with the literature as regards the $p-$dependent convergence observed for IgA discretizations combined with multigrid. Generally speaking, more smoothing steps and/or intricate smoothers are needed in order to counteract the increasing number of iterations for higher-order IgA schemes. \\
% Both $DC_{MG}^{1}$ and $D_{\varepsilon}C_{MG}^{1}$ result in a dependency of the iteration numbers on $p$. {\color{red} These results are in line with the spectral analysis from Section \ref{spec}, in particular Figure \ref{fig:p} and Figure \ref{fig:p}, }However, with respect to $DC_{MG}^{1}$, starting from $k = 10^5$ the dependency becomes inversely related and we observe that the number of iterations starts increasing with $k$ but decreasing with $p$. Once we add the weight parameter $\varepsilon$ to the deflation preconditioner we obtain $k$-independent convergence up to $10^6$. \\
% The dependency on the approximation order $p$ can be further mitigated by applying more multigrid cycles. In fact, application of $10$ multigrid cycles leads to $p$-independent convergence and shows identical results to inverting CSLP exactly; see Table \ref{tab:2}. Note, however, that application of $D_{\varepsilon}C_{MG}^{10}$ is more expensive compared to the application of $D_{\varepsilon}C_{MG}^{1}$ as we use more V-cycles in order to obtain a fair approximation of the CSLP. This result, however, is in line with the literature as regards the $p-$dependent convergence observed for IgA discretizations combined with multigrid. Generally speaking, more smoothing steps and/or intricate smoothers are needed in order to counteract the increasing number of iterations for higher-order IgA schemes. \\
\begin{table}[ht!]
\centering
\scalebox{0.75}{
\begin{tabular}{c|ccc|ccc|ccc|ccc|ccc|}
& \multicolumn{3}{c|}{$k=10^2$} &  \multicolumn{3}{c|}{$k=10^3$} & \multicolumn{3}{c|}{$k=10^4$} &  \multicolumn{3}{c|}{$k=10^5$}  &  \multicolumn{3}{c|}{$k=10^6$} \\ \hline
& \multicolumn{3}{c|}{$N = 161$} &  \multicolumn{3}{c|}{$N = 1601$} & \multicolumn{3}{c|}{$N = 16001$} &  \multicolumn{3}{c|}{$N = 160001$}  &  \multicolumn{3}{c|}{$N = 1600001$} \\
     & $DC_{MG}^{1}$ & $D_{\varepsilon}C_{MG}^{1}$  & $D_{\varepsilon}C_{MG}^{10}$    & $DC_{MG}^{1}$ & $D_{\varepsilon}C_{MG}^{1}$  & $D_{\varepsilon}C_{MG}^{10}$    & $DC_{MG}^{1}$   & $D_{\varepsilon}C_{MG}^{1}$  & $D_{\varepsilon}C_{MG}^{10}$   & $DC_{MG}^{1}$ &  $D_{\varepsilon}C_{MG}^{1}$   & $D_{\varepsilon}C_{MG}^{10}$  & $DC_{MG}^{1}$ & $D_{\varepsilon}C_{MG}^{1}$  & $D_{\varepsilon}C_{MG}^{10}$  \\ \hline  
$p=1$ & $7$	 &$7$ &$5$  & $7$  &$7$ & $5$	  &$13$ & $7$	  &$5$ & $50$	  &$7$  &$5$ &$*$ &$10$ &$5$ \\
$p=2$ & $5$	 &$5$ &$5$  & $6$  &$5$ & $5$	  &$10$ & $5$	  &$5$ & $28$	  &$5$  &$5$ &$*$ &$5$ &$5$\\
$p=3$ & $6$	 &$6$ &$5$  & $6$  &$6$ & $5$	  &$8$ & $6$	  &$5$ & $22$	  &$6$  &$5$ &$*$ &$6$ &$5$\\
$p=4$ & $9$  &$9$ &$5$  & $9$  &$9$ & $5$	  &$10$ & $10$	  &$5$ & $19$	  &$9$  &$5$ &$74$ &$9$ &$5$\\
$p=5$ & $16$ &$16$ &$5$ & $16$ &$16$ & $5$	  &$13$ & $15$	  &$5$ & $21$	  &$15$  &$5$ &$46$ &$15$ &$5$\\
\end{tabular}
}
\caption{Number of (preconditioned) GMRES iterations to reach convergence for MP 1-B. Here we combine the two-level deflation (D) using \textbf{quadratic} Bezier curves with the CSLP. The shift $\beta_2$ has been set to 1. CSLP has been inverted using $C_{MG}^1$ and $C_{MG}^{10}$ respectively.}
\label{tab:4}
\end{table}

As mentioned previously, for a finite difference scheme, it has been shown that the deflation preconditioner without CSLP could also lead to close to wave number independent convergence. Thus, analogously, we perform a similar test to examine how well the deflation preconditioner performs with no other preconditioner. We will distinguish two cases; ADP without weight parameter $D$ and ADP with weight parameter $D_{\varepsilon}$. For $D$, the results are reported in Table \ref{tab:2}, where we compare the number of iterations to the number of iterations obtained by using the (exactly inverted) CSLP with shift $k^{-1}$ $(C_{ex})$. Note that, the exactly inverted CSLP leads to iteration numbers independent of both $k$ and $p$. In absence of the weight parameter, the number of GMRES iterations preconditioned with $D$ increases with $k$ and $p$ for wave numbers $k < 10^5$. These results are similar to the ones reported in Table \ref{tab:4}, where we observed a similar effect for $DC_{MG}^{1}$. The observed number of iterations is also in agreement with the spectral analysis from Fig \ref{fig:pa2} in Section \ref{spec}. It has been shown that as the accuracy of ADP decreases, the projection error increases, and the eigenvalues are not accurately projected onto the origin. As a result, the number of iterations is expected to increase with $k$. However, we did observe that this effect is less pronounced for larger values of $p$, which is why we obtain better convergence for larger values of $k$ when $p \geq 4$. \\
\begin{table}[ht!]
\centering
\scalebox{1.00}{
\begin{tabular}{c|cc|cc|cc|cc|cc|}
& \multicolumn{2}{c|}{$k=10^2$} &  \multicolumn{2}{c|}{$k=10^3$} & \multicolumn{2}{c|}{$k=10^4$} &  \multicolumn{2}{c|}{$k=10^5$}  &  \multicolumn{2}{c|}{$k=10^6$} \\ \hline
& \multicolumn{2}{c|}{$N = 161$} &  \multicolumn{2}{c|}{$N = 1601$} & \multicolumn{2}{c|}{$N = 16001$} &  \multicolumn{2}{c|}{$N = 160001$}  &  \multicolumn{2}{c|}{$N = 1600001$} \\
     & $D$   & $C_{ex}$    & $D$   & $C_{ex}$    & $D$   & $C_{ex}$  & $D$   & $C_{ex}$ &  $D$   & $C_{ex}$    \\ \hline  
$p=1$ & $9$	  &$5$ & $8$	  &$5$ & $13$	  &$5$ & $49$	  &$5$ & $*$	  &$5$   \\
$p=2$ & $7$	  &$5$ & $6$	  &$5$ & $10$	  &$5$ & $28$	  &$5$ & $*$	  &$5$   \\
$p=3$ & $8$	  &$5$ & $8$	  &$5$ & $10$	  &$5$ & $20$	  &$5$ & $*$	  &$5$   \\
$p=4$ & $13$  &$5$ & $13$  &$5$ & $13$	  &$5$ & $20$	  &$5$ & $68$	  &$5$   \\
$p=5$ & $19$  &$5$ & $19$  &$5$ & $16$	  &$5$ & $25$	  &$5$ & $48$	  &$5$   \\
\end{tabular}
}
\caption{Number of (preconditioned) GMRES iterations to reach convergence for MP 1-B. Here we use GMRES with either two-level deflation (D) using \textbf{quadratic} Bezier curves or exact inverse of CSLP $C_{ex}$ using $\beta_2 = k^{-1}$. * indicates that the number of max iterations (100) has been reached without convergence.}
\label{tab:2}
\end{table}
Table \ref{tab:3} contains the same comparison, however we use the deflation preconditioner $D_{\varepsilon}$. We report the number of (preconditioned) GMRES for both preconditioners. Note that, the exactly inverted CSLP leads to iteration numbers independent of both $k$ and $p$. In absence of the weight parameter, the number of GMRES iterations adopting the deflation preconditioner increases with $k$ and decreases with $p$ starting from $k = 10^5$. These results are similar to the ones reported in Table \ref{tab:4}, where we observed a similar effect for $DC_{MG}^{1}$. Adding the weight parameter significantly improves the convergence of the GMRES method with respect to $k-$dependent convergence. In particular, wave number independent convergence is observed for values of $k$ up to $10^6$. This is in line with the spectral analysis from Fig \ref{fig:pa1} in Section \ref{spec}. There, we observed that an accurate interpolation scheme ensures that half of the eigenvalues are mapped onto the origin and the spectrum remains as clustered as possible. However, for $p=5$ we observed that the smallest and largest eigenvalue lie further away, which could explain the $p-$dependent convergence, and in particular the higher number of iterations observed for $p = 5$. Thus, similar to multigrid solvers, deflation based solvers are also subjected to $p-$dependent convergence. The effect can be circumvented by combining both methods and increasing the number of V-cycles. 
\begin{table}[ht!]
\centering
\scalebox{1.00}{
\begin{tabular}{c|cc|cc|cc|cc|cc|}
& \multicolumn{2}{c|}{$k=10^2$} &  \multicolumn{2}{c|}{$k=10^3$} & \multicolumn{2}{c|}{$k=10^4$} &  \multicolumn{2}{c|}{$k=10^5$}  &  \multicolumn{2}{c|}{$k=10^6$} \\ \hline
& \multicolumn{2}{c|}{$N = 161$} &  \multicolumn{2}{c|}{$N = 1601$} & \multicolumn{2}{c|}{$N = 16001$} &  \multicolumn{2}{c|}{$N = 160001$}  &  \multicolumn{2}{c|}{$N = 1600001$} \\
     & $D_{\varepsilon}$   & $C_{ex}$    & $D_{\varepsilon}$   & $C_{ex}$    & $D_{\varepsilon}$   & $C_{ex}$  & $D_{\varepsilon}$   & $C_{ex}$ &  $D_{\varepsilon}$   & $C_{ex}$    \\ \hline  
$p=1$ & $9$	  &$5$ & $9$	  &$5$ & $9$	  &$5$ & $9$	  &$5$ & $11$	  &$5$   \\
$p=2$ & $5$	  &$5$ & $5$	  &$5$ & $5$	  &$5$ & $5$	  &$5$ & $5$	  &$5$   \\
$p=3$ & $8$	  &$5$ & $8$	  &$5$ & $8$	  &$5$ & $8$	  &$5$ & $8$	  &$5$   \\
$p=4$ & $13$  &$5$ & $13$  &$5$ & $13$	  &$5$ & $11$	  &$5$ & $13$	  &$5$   \\
$p=5$ & $20$  &$5$ & $20$  &$5$ & $20$	  &$5$ & $19$	  &$5$ & $20$	  &$5$   \\
\end{tabular}
}
\caption{Number of (preconditioned) GMRES iterations to reach convergence for MP 1-B. Here we use GMRES with two-level deflation ($D_{\varepsilon}$) using \textbf{quadratic} Bezier curves. $C_{ex}$ uses the shift $\beta_2 = k^{-1}$ and is inverted exactly. * indicates that the number of max iterations (100) has been reached without convergence.}
\label{tab:3}
\end{table}
\FloatBarrier

\subsubsection{Two-dimensional model problems}
\paragraph{MP 2-A}
\FloatBarrier
In the previous subsection, it was observed that combining the deflation preconditioner $D_{\varepsilon}$ with the approximated CSLP $C^{j}_{MG}$ yields the best results in terms of iteration numbers. In this subsection, we apply this preconditioner to MP 2-A, the natural extension of MP 1-B to two dimensions. In particular, CPU timings are determined to obtain a fair comparison in terms of computational costs. \\ \\
Table \ref{tab:5} compares $DC_{MG}^1$ and $D_{\varepsilon}C_{MG}^{12}$ with the exactly inverted CSLP $C_{ex}$. For $D_{\varepsilon}C_{MG}^{12}$, we obtain close to $k$- and $p$- independent convergence. Only for $p=5$, the number of iterations increases. Here, $3$ pre- and post-smoothing steps and a shift of $\beta_2=4.2$ are adopted. For the $C_{ex}$ preconditioner, a shift of $(3k)^{-1}$ has been adopted. Both the shift $k^{-1}$ as well the shift $\beta_2 = (3k)^{-1}$ does not lead to wave number independent convergence. In fact, $C_{ex}$
 uses more iterations for $p < 5$ in most cases. This can be explained by the fact that we are using Dirichlet boundary conditions, which are known to cause a less favorable distribution of the eigenvalues compared to the use of Sommerfeld radiation conditions \cite{van2007spectral}. In particular, keeping the shift $k^{-2}$ results in wave number independent convergence but leads to very uneconomical systems, which are close to the original coefficient matrix. 
\begin{table}[ht!]
\centering
\scalebox{0.75}{
\begin{tabular}{c|ccc|ccc|ccc|ccc|ccc|}
& \multicolumn{3}{c|}{$k=50$} &  \multicolumn{3}{c|}{$k=100$} & \multicolumn{3}{c|}{$k=150$} &  \multicolumn{3}{c|}{$k=200$}  &  \multicolumn{3}{c|}{$k=250$} \\ \hline
& \multicolumn{3}{c|}{$N = 6241$} &  \multicolumn{3}{c|}{$N = 25281$} & \multicolumn{3}{c|}{$N = 57121$} &  \multicolumn{3}{c|}{$N = 101761$}  &  \multicolumn{3}{c|}{$N = 159201$} \\
     & $DC_{MG}^{1}$ & $D_{\varepsilon}C_{MG}^{12}$  & $C_{ex}$    & $DC_{MG}^{1}$ & $D_{\varepsilon}C_{MG}^{12}$  & $C_{ex}$    & $DC_{MG}^{1}$   & $D_{\varepsilon}C_{MG}^{12}$  & $C_{ex}$   & $DC_{MG}^{1}$ &  $D_{\varepsilon}C_{MG}^{12}$   & $C_{ex}$  & $DC_{MG}^{1}$ & $D_{\varepsilon}C_{MG}^{12}$  & $C_{ex}$  \\ \hline  
$p=1$ & $7$	 &$7$ &$7$  & $8$  &$7$ & $8$	  &$12$ & $12$	  &$10$ & $8$	  &$8$  &$9$ &$12$ &$9$ &$10$ \\
$p=2$ & $10$	 &$7$ &$7$  & $10$  &$7$ & $8$	  &$10$ & $7$	  &$8$ & $11$	  &$8$  &$11$ &$12$ &$8$ &$10$\\
$p=3$ & $18$	 &$6$ &$6$  & $20$  &$9$ & $8$	  &$18$ & $7$	  &$7$ & $20$	  &$7$  &$11$ &$19$ &$7$ &$10$\\
$p=4$ & $36$  &$7$ &$6$  & $36$  &$7$ & $8$	  &$36$ & $7$	  &$7$ & $36$	  &$7$  &$11$ &$37$ &$7$ &$10$\\
$p=5$ & $85$ &$20$ &$7$ & $86$ &$21$ & $8$	  &$87$ & $21$	  &$7$ & $86$	  &$21$  &$11$ &$21$ &$21$ &$10$\\
\end{tabular}
}
\caption{Number of (preconditioned) GMRES iterations to reach convergence for MP 2-A. Here we combine the two-level deflation (D) using \textbf{quadratic} Bezier curves with CSLP. CSLP has been inverted using $C_{MG}^1$ and $C_{MG}^{12}$ respectively where the shift has been set to $\beta_2 = 1$ and $\beta_2 = 4.2$ respectively. When using $C_{ex}$, the shift has been set to $\beta_2 = 3k^{-1}$.}
\label{tab:5}
\end{table}

\begin{figure}[ht!]
\centering
   \begin{tikzpicture}
    \begin{axis}[width=\textwidth,
       height=0.45\textwidth,nodes near coords={\pgfmathprintnumber[fixed,fixed zerofill,precision=1]{\pgfplotspointmeta}},point meta=rawy,ybar,every node near coord/.append style={rotate=90, anchor=west,font=\small},legend columns=2,legend pos=north west, xmin= 0,xmax=62,ymode=log,ymin=1,ymax=2000,ylabel=CPU Time (s),line width=0.25pt,bar width=0.2cm,bar shift=8mm,xticklabels=\emptynodes]
       \addplot[fill=blue!20] plot coordinates 
      {
        (1,    1.70)
        (13,   7.52)
        (25,   14.15)
        (37,   31.39)
        (49,   81.96)
        };
      \addlegendentry{$p=2, DC$}
       \addplot[fill=blue!40] plot coordinates 
      {
        (2,    1.17)
        (14,   8.18)
        (26,   15.41)
        (38,   43.21)
        (50,   133.81)
        };
      \addlegendentry{$p=2, C$}
        \addplot[fill=orange!20] plot coordinates 
      {
        (3,    2.15)
        (15,   8.94)
        (27,   18.22)
        (39,   35.16)
        (51,   62.98)
        };
      \addlegendentry{$p=3, DC$}
      \addplot[fill=orange!40] plot coordinates 
      {
        (4,    2.15)
        (16,   13.45)
        (28,   32.15)
        (40,   94.48)
        (52,   161.16)
        };
      \addlegendentry{$p=3, C$}
       \addplot[fill=green!20] plot coordinates 
      {
        (5,    3.43)
        (17,   13.73)
        (29,   28.33)
        (41,   55.43)
        (53,   95.03)
        };
      \addlegendentry{$p=4, DC$}       
      \addplot[fill=green!40] plot coordinates 
      {
        (6,    2.93)
        (18,   20.18)
        (30,   58.81)
        (42,   166.31)
        (54,   295.40)
        };
      \addlegendentry{$p=4, C$}
      \addplot[fill=red!20] plot coordinates 
      {
        (7,    7.36)
        (19,   33.97)
        (31,   88.24)
        (43,   156.23)
        (55,   264.44)
        };
      \addlegendentry{$p=5, DC$}
      \addplot[fill=red!40] plot coordinates 
      {
        (8,    4.88)
        (20,   33.29)
        (32,   126.79)
        (44,   302.74)
        (56,   933.16)
        };
      \addlegendentry{$p=5, C$}
      \end{axis}
       \draw (1.95,-0.45)  node{$k=50$};
        \draw (5.10,-0.45)  node{$k=100$};
        \draw (8.25,-0.45) node{$k=150$};
        \draw (11.40,-0.45) node{$k=200$};
        \draw (14.55,-0.45) node{$k=250$};
\end{tikzpicture}
\caption{CPU-time in seconds (s) for $p=2$ to $p=5$ for MP 2-A. The plot contains the timings for $k=50,100,150,200$ and $k=250$. DC stands for $D_{\varepsilon}C_{MG}^{12}$ and $C$ stands for $C_{ex}$ using $\beta_2 = (3k)^{-1}$.}
\label{fig:1b}
\end{figure}
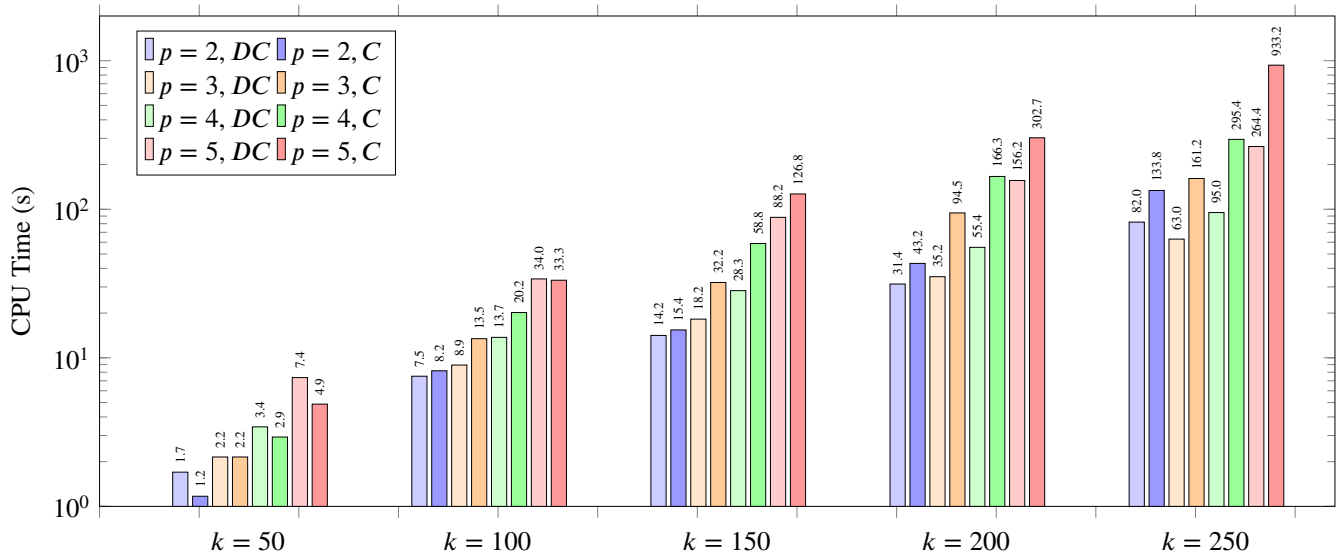

Figure \ref{fig:1b} shows the corresponding CPU times to reach convergence with the GMRES method when applying $D_{\varepsilon}C_{MG}^{12}$ and $C_{ex}$ as a preconditioner. The CPU-timings have been obtained using the Matlab 2019b 'tic toc' command. For $k=50$, inverting the CSLP preconditioner exactly leads to the lowest CPU times for all values of $p$ considered. However, from $k=150$ already, the opposite holds: $D_{\varepsilon}C_{MG}^{12}$ is computationally more efficient compared to the exact CSLP preconditioner. This effect becomes more pronounced as $k$ increases. Thus, the larger $k$, the larger the computational speedup of the deflated preconditioned solver relative to the solver using the exact inversion of the CSLP combined with a small complex shift. 

\FloatBarrier

\begin{figure}[ht!]
    \centering
    \includegraphics[scale=0.41]{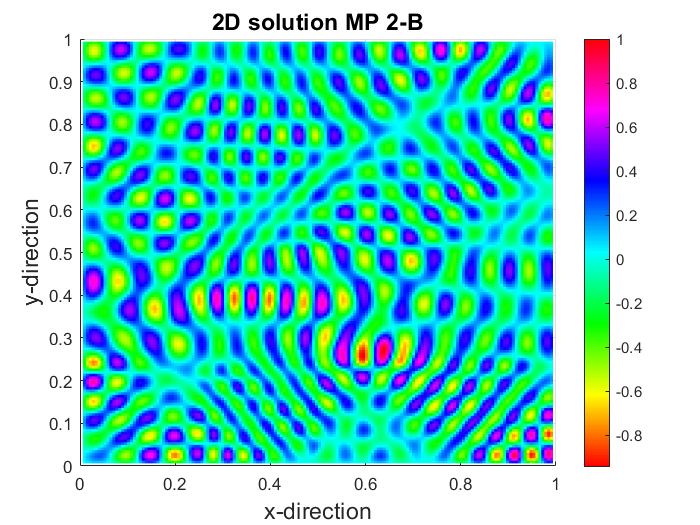}
    \includegraphics[scale=0.41]{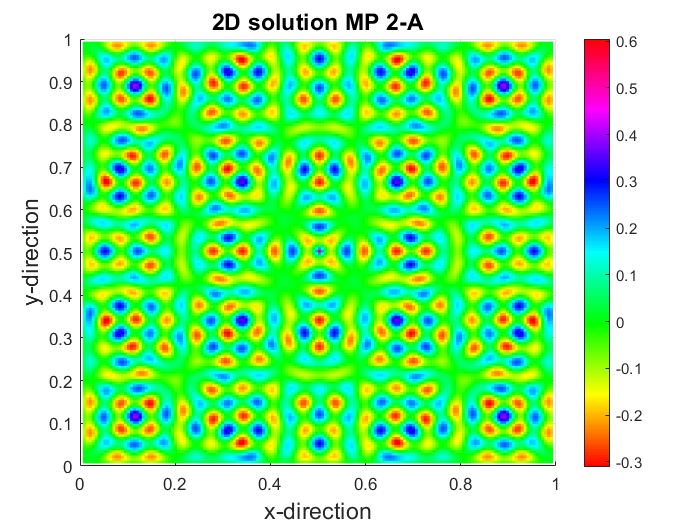}
    \caption{Real part of the two-dimensional numerical solution for the non-constant wave number $k(x,y)$ where $k = 100$ (left) and $k = 100$ (right).}
    \label{fig:solk100}
\end{figure}

\paragraph{MP 2-B}
\FloatBarrier
Finally, we consider model problem MP 2-B, where the wave number is non-constant and given by a two-dimensional step function. This is an important benchmark as some solvers only perform successfully when a constant wave number is used. Moreover, it allows for testing whether the numerical solver can deal with sharp disruptions in the underlying velocity, which is the main focus of this section. 
In Figure \ref{fig:solk100} we have plotted the variable (left) and constant (right) solution for MP 2-A and MP 2-B respectively using $k = 100$ as a base wave number. The step-function used to vary $k$ throughout the numerical domain is observed to disrupt the symmetric pattern observed for $k = 100$ (right). 

Table \ref{tab:nonconsk} shows the number of GMRES iterations needed to reach convergence when $DC_{MG}^{12}$ and $C_{ex}$ are applied as a preconditioner. With respect to $p-$dependent convergence, the number of iterations slightly varies with $p$ for both preconditioned systems. In contrast to MP 2-A, however, we also observe a small increase in the number of iterations as $k$ increases for both preconditioned systems. However, in terms of iterations, the deflated preconditioned system needs less iterations compared to the system using the exact inversion of the CSLP and a very small complex shift. Unlike the results from the constant wave number model problem, we therefore report weakly dependent convergence on $k$. However, note that for $p = 5$, the convergence appears to resemble wave number independent convergence. We do note that using the deflation preconditioner combined with the multigrid approximation of the CSLP, the number of iterations could be improved by using more V-cycles. These are relatively cheap in terms of computational costs as they are of order $\mathcal{O}(N)$ FLOPs and given that the diagonal scaled Jacobi smoother is used. 
\begin{table}[ht!]
\centering
\scalebox{1.00}{
\begin{tabular}{c|cc|cc|cc|cc|cc|}
& \multicolumn{2}{c|}{$k=50$} &  \multicolumn{2}{c|}{$k=100$} & \multicolumn{2}{c|}{$k=150$} &  \multicolumn{2}{c|}{$k=200$}  &  \multicolumn{2}{c|}{$k=250$} \\ \hline
& \multicolumn{2}{c|}{$N = 6241$} &  \multicolumn{2}{c|}{$N = 25281$} & \multicolumn{2}{c|}{$N = 57121$} &  \multicolumn{2}{c|}{$N = 101761$} &  \multicolumn{2}{c|}{$N = 159201$} \\ 
& $DC_{MG}^{12}$   & $C_{ex}$    & $DC_{MG}^{12}$   & $C_{ex}$    & $DC_{MG}^{12}$   & $C_{ex}$  & $DC_{MG}^{12}$   & $C_{ex}$ &  $DC_{MG}^{12}$   & $C_{ex}$    \\ \hline  
$p=1$ & $13$  &$12$ & $16$	  &$19$ & $22$	  &$24$ & $25$	  &$27$ & $29$	  &$28$   \\
$p=2$ & $13$  &$13$ & $16$	  &$20$ & $20$	  &$24$ & $25$	  &$29$ & $32$	  &$36$   \\
$p=3$ & $10$  &$13$ & $11$	  &$16$ & $14$	  &$23$ & $15$	  &$28$ & $20$	  &$39$   \\
$p=4$ & $10$  &$13$ & $13$    &$20$ & $12$	  &$22$ & $13$	  &$26$ & $19$	  &$38$   \\
$p=5$ & $18$  &$13$ & $19$    &$16$ & $17$	  &$23$ & $21$	  &$29$ & $20$	  &$39$   \\
\end{tabular}
}
\caption{Number of (preconditioned) GMRES iterations to reach convergence for MP 2-B. Here we combine two-level deflation using \textbf{quadratic} Bezier curves with CSLP $(DC_{MG}^{12})$. For $p < 5$ we use 3 pre- and post smoothing steps, whereas for $p=5$ we use 2 pre- and post smoothing steps. CSLP has been inverted using $C_{MG}^{12}$ where the shift has been set to $\beta_2 = 4.2$. When using $C_{ex}$, the shift has been set to $\beta_2 = {(3k)}^{-1}$ and CSLP is inverted exactly.}
\label{tab:nonconsk}
\end{table}

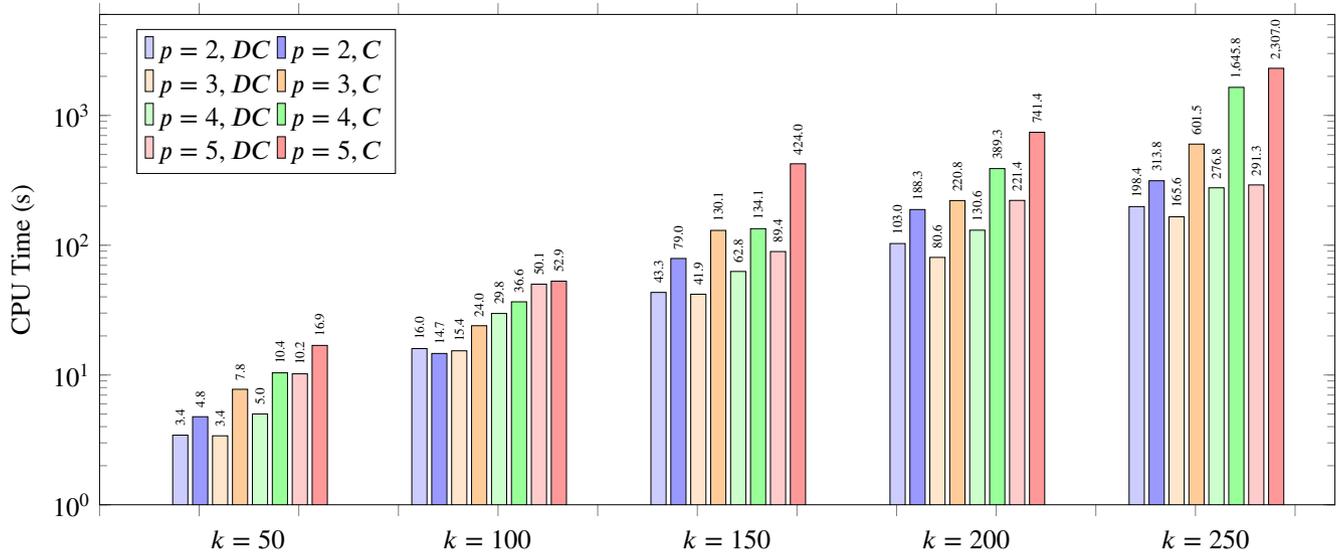
\begin{figure}[ht!]
\centering
   \begin{tikzpicture}
    \begin{axis}[width=\textwidth,
       height=0.45\textwidth,nodes near coords={\pgfmathprintnumber[fixed,fixed zerofill,precision=1]{\pgfplotspointmeta}},point meta=rawy,ybar,every node near coord/.append style={rotate=90, anchor=west,font=\small},legend columns=2,legend pos=north west, xmin= 0,xmax=62,ymode=log,ymin=1,ymax=6000,ylabel=CPU Time (s),line width=0.25pt,bar width=0.2cm,bar shift=8mm,xticklabels=\emptynodes]
       \addplot[fill=blue!20] plot coordinates 
      {
        (1,    3.44)
        (13,   15.98)
        (25,   43.29)
        (37,   102.99)
        (49,   198.43)
        };:(
        
      \addlegendentry{$p=2, DC$}
       \addplot[fill=blue!40] plot coordinates 
      {
        (2,    4.77)
        (14,   14.65)
        (26,   79.04)
        (38,   188.32)
        (50,   313.83)
        };
      \addlegendentry{$p=2, C$}
        \addplot[fill=orange!20] plot coordinates 
      {
        (3,    3.40)
        (15,   15.36)
        (27,   41.92)
        (39,   80.55)
        (51,   165.55)
        };
      \addlegendentry{$p=3, DC$}
      \addplot[fill=orange!40] plot coordinates 
      {
        (4,    7.76)
        (16,   23.97)
        (28,   130.05)
        (40,   220.83)
        (52,   601.52)
        };
      \addlegendentry{$p=3, C$}
       \addplot[fill=green!20] plot coordinates 
      {
        (5,    5.01)
        (17,   29.79)
        (29,   62.78)
        (41,   130.55)
        (53,   276.83)
        };
      \addlegendentry{$p=4, DC$}       
      \addplot[fill=green!40] plot coordinates 
      {
        (6,    10.39)
        (18,   36.64)
        (30,   134.08)
        (42,   389.32)
        (54,   1645.81)
        };
      \addlegendentry{$p=4, C$}
      \addplot[fill=red!20] plot coordinates 
      {
        (7,    10.21)
        (19,   50.11)
        (31,   89.36)
        (43,   221.40)
        (55,   291.32)
        };
      \addlegendentry{$p=5, DC$}
      \addplot[fill=red!40] plot coordinates 
      {
        (8,    16.90)
        (20,   52.87)
        (32,   423.97)
        (44,   741.36)
        (56,   2306.96)
        };
      \addlegendentry{$p=5, C$}
      \end{axis}
       \draw (1.95,-0.45)  node{$k=50$};
        \draw (5.10,-0.45)  node{$k=100$};
        \draw (8.25,-0.45) node{$k=150$};
        \draw (11.40,-0.45) node{$k=200$};
        \draw (14.55,-0.45) node{$k=250$};
\end{tikzpicture}
\caption{CPU-time in seconds (s) for $p=2$ to $p=5$ for MP 2-B. The plot contains the timings for $k=50,100,150,200$ and $k=250$. DC stands for $D_{\varepsilon}C_{MG}^{12}$ and $C$ stands for $C_{ex}$ using $\beta_2 = (3k)^{-1}$.}
\label{fig:k1}
\end{figure}

The corresponding CPU timings are provided in Figure \ref{fig:k1}. The combination of deflation and the approximated deflation preconditioner $(DC_{MG}^{12})$ is cheaper for all values of $p$ and $k$. Hence, already for moderate values of $k$, applying the CSLP preconditioner exactly is more expensive. Note that, for higher values of $k$, the difference between both approaches also becomes more visible in terms of CPU timings. This effect will only be magnified in 3D-applications.
%we take p = 2 till 5 for the graph
\FloatBarrier

\section{Conclusion}\label{concl}
In this work, we focus on the combination of IgA discretized linear systems with a state-of-the-art iterative solver using deflation and a geometric multigrid method. In particular, we extend the line of research set out by \cite{diwan2020iterative}, where it was shown that the use of IgA reduces the pollution error significantly compared to $p-$order FEM. The authors have shown that the use of the exact inverse of the CSLP preconditioner with a small complex shift, yields wave number independent convergence for moderate values of $k$. Instead of inverting the CSLP exactly and using a small complex shift, we use a standard multigrid method to approximate its inverse and combine it with a two-level deflation preconditioner to accelerate the convergence of GMRES.
We use a large complex shift in order to ensure that the multigrid algorithm does not diverge.

% {\color{red}The pollution error causes the numerical solution to deteriorate as the wave number $k$ increases and can not be avoided. As a result, a growing branch of research has focused on improving discretization techniques, such as using IgA discretizations for the Helmholtz equation. 
%  When it comes to solving the Helmholtz equation, iterative solution methods are preferred as the resulting linear systems are often large, indefinite and complex. However, the number of iterations increases with the wave number $k$ and the search for scalable solvers in terms of iterations and CPU-timings remains an active research topic. {\color{cyan}{Dit voelt alemaal een herhaling. Misschien niet in de conclusie ook nog? Wel: kort nog beschrijven wat ze in [11] hebben gedaan (naast de pollution error)}}}

% The use of Dirichlet boundary conditions has enabled us to study the 'worse-case' scenario for the IgA discretizations as the natural damping effect from, for example the Sommerfeld boundary conditions, is now omitted from the spectrum. {\color{cyan}{Misschien te uitgebreid voor een conclusie?}} 
%This accelerates the convergence of GMRES as this Krylov subspace method is more susceptible to having small and scattered eigenvalues.
The use of deflation techniques is motivated by studying the spectrum of the preconditioned systems. Deflation projects the unwanted negative and near-zero eigenvalues corresponding to the smooth eigenmodes onto zero, thereby accelerating the convergence of GMRES. Our spectral analysis shows that for increasing $k$ and $p$, the spectrum remains well-clustered. This is supported by the numerical results in 1D as the number of iterations remains $k$- and $p$-independent for $kh$ constant. If we exclude the CSLP, we obtain $k$ independent convergence and the number of iterations increases slightly with $p$. 

When deflation is combined with CSLP, the number of iterations weakly depends on $k$ and $p$ for $kh$ constant in the 2D case. Starting from $k = 150$, the deflation based preconditioner combined with the approximate inverse of the CSLP outperforms the exact inversion of the CSLP with shift $\beta_2 = (3k)^{-1}$ in terms of CPU-timings. The obtained speed-up becomes more significant as the wave number $k$ increases. Results for the highly varying non-constant wave number model show a slight dependence on $k$ but an inversely related dependence on $p$ as the wave number increases. Even for this model problem, the proposed solver outperforms in terms of number of iterations and CPU-timings, when compared to the use of the exact inversion of the CSLP with a small complex shift. 
%\section*{Acknowledgments}

%\subsection*{Author contributions}

%\subsection*{Financial disclosure}
%None reported.

%\subsection*{Conflict of interest}
%The authors declare no potential conflict of interests.
\bibliography{HelmholtzIgA}%

\end{document}